\PassOptionsToPackage{usenames,dvipsnames}{xcolor}
\PassOptionsToPackage{verbose=false}{myPreamble}
\documentclass[10pt, USenglish, reqno]{article}

\makeatletter
	\PassOptionsToPackage{nameinlink,capitalize}{cleveref}
	\PassOptionsToPackage{breakspaces}{clevethm}
	\PassOptionsToPackage{inline}{enumitem}
	\PassOptionsToPackage{bookmarksdepth=4}{hyperref}
	\usepackage[bottom]{footmisc}
	\usepackage{tabularx}
	\usepackage[most]{tcolorbox}
	\usepackage[
		arxiv,
		alglinenumber,
		excludethm=ass,
		shortthmnames,
	]{myPreamble}

	\pgfplotsset{compat=1.16}
	\usepgfplotslibrary{fillbetween}

	\crefname{footnote}{footnote}{footnotes}

	\Alignfalse%
	\newcommand{\NL}[1][]{\ifstrempty{#1}{,}{#1}\ }
	\newcommand{\AND}{ \(\cdot\) }
	\newcommand{\subclass}[1]{\par\noindent{\def\and{\unskip,\ }{\bf AMS subject classifications. }#1}\par}
	\let\oldand\and
	\def\and{\texorpdfstring{\oldand}{, }}%

	\AtBeginDocument{%
		\let\oldmybox\mybox
		\let\oldendmybox\endmybox
		\@namedef{mybox*}{\oldmybox}
		\@namedef{endmybox*}{\oldendmybox}
	}

	\Crefname{figure}{Figure}{Figures}

	\theoremstyle{plain}
	\CreateNewTheorem{ass}{Assumption}
	\newlist{claims}{enumerate}{1}
	\setlist[claims]{%
		label={\it Claim \alph*.},
		ref=\thedummythm.\alph*,
		wide=0pt,
		widest=99,
		leftmargin=\parindent,
	}
	\setlistlabel[enumerate,1]{{%
		\arabic*%
	}}
	\setlistlabel[enumerateq,1]{{%
		\alph*%
	}}
	\setlistlabel[enumerator,1]{{%
		\Alph.%
	}}
	\setlistlabel[enumeratass,1]{{%
		{\sc a}\oldstylenums{\arabic*}%
	}}
	\setlistlabel[enumeratprop,1]{{%
		{\scriptsize p}\oldstylenums{\arabic*}%
	}}
	\def\refenumeratelabeli{\arabic*}
	\def\refenumerateqlabeli{\alph*}
	\def\refenumeratorlabeli{\Alph*}
	\def\refenumeratasslabeli{\textsc{a}\oldstylenums{\arabic*}}
	\def\refenumeratproplabeli{\textsc{p}\oldstylenums{\arabic*}}
	\def\refenumeratePrefixi{.}
	\def\refenumerateqPrefixi{.}
	\def\refenumeratorPrefixi{.}
	\def\refenumeratassPrefixi{.}
	\def\refenumeratpropPrefixi{.}
	\def\enumerateSuffixi{.}
	\def\enumerateqSuffixi{.}
	\def\enumeratorSuffixi{.}
	\def\enumeratassSuffixi{.}
	\def\enumeratpropSuffixi{.}

	\crefname{claimsi}{claim}{claims}
	\crefformat{claimsi}{claim #2#1#3}

	\newlist{algsteps}{enumerate}{1}
	\setlist[algsteps,1]{%
		label=\oldstylenums{\arabic*}:,
		ref=\oldstylenums{\arabic*},
		topsep=0pt,
	}
	\crefname{algstepsi}{step}{steps}

	\renewcommand{\algfont}[1]{\textsc{#1}}
	\crefname{ALG@line}{step}{steps}
	
	\renewcommand{\theALG@line}{{\bf\thealgorithm}.{\oldstylenums{\arabic{ALG@line}}}}
	\newlist{conditions}{enumerate}{1}
	\setlist[conditions]{%
		label=(\alph*),
		ref=\theALG@line(\alph*),
		leftmargin=*,
		nosep,
		noitemsep,
	}
	\crefname{conditionsi}{condition}{conditions}
	\crefformat{conditionsi}{condition #2#1#3}

	\let\olditem\item
	\newcommand{\@itemwitharg}[1]{\ref{#1})~}%
	\AtBeginEnvironment{proofitemize}{%
		\renewcommand{\item}{%
			\olditem
			\@ifnextchar\bgroup{\@itemwitharg}{}%
		}%
	}

	\newtcolorbox{mybox}[1][]{%
		left=0pt,
		right=0pt,
		top=0pt,
		bottom=0pt,
		colback=MidnightBlue!5,
		colframe=MidnightBlue!10,
		width=\dimexpr\textwidth\relax,
		enlarge left by=0mm,
		boxsep=3pt,
		arc=5pt,outer arc=5pt,
		after=\ignorespacesafterend\par\noindent,
		#1
	}

	\newcommand{\ippanoc}{\@ifnextchar({\@ippanoc@arg}{\@ippanoc@noarg}}
	\def\@ippanoc@noarg{%
		\text{\texorpdfstring{\hyperref[alg:PANOC+]{IP-PANOC\(^+\)}}{PANOC+}}%
	}
	\def\@ippanoc@arg(#1){%
		\@ippanoc@noarg\texorpdfstring{\ensuremath{(#1)}}{}%
	}
	\newcommand{\ipfb}{\text{\texorpdfstring{\hyperref[alg:FB]{IP-FB}}{IP-FB}}}

	\newcommand{\x}{z}%

	\newcommand{\Res}{\@ifnextchar_\Res@sub\Res@no@sub}
	\def\Res@sub_#1{\operatorname{R}_{#1}}
	\def\Res@no@sub{\operatorname{R}_{\FBstepsize}}

	\newcommand{\T}{%
		\def\@command{%
			\ifx\@subarg\relax
				\operatorname{T}^{\text{\sc fb}}_{\IPmu,\FBstepsize}
			\else
				\operatorname{T}^{\text{\sc fb}}_{\@subarg}
			\fi
		}%
		\@getscripts%
	}

	\newcommand{\iter}{k}

	\renewcommand{\FBstepsize}{\@ifstar\@@FBstepsize\@FBstepsize}
	\newcommand{\@FBstepsize}{\gamma_{\iter}}
	\newcommand{\@@FBstepsize}{\gamma_{\iter-1}}
	\renewcommand{\FBsmooth}{f_{\IPmu}}

	\def\IPmu{\mu}
	\newcommand{\FBE}{\@ifstar\@@FBE\@FBE}
	\newcommand{\@FBE}{%
		\ifx\IPmu\empty
			\varphi_{\FBstepsize}^{\text{\sc fb}}
		\else
			\varphi_{\IPmu,\FBstepsize}^{\text{\sc fb}}
		\fi
	}
	\newcommand{\@@FBE}{%
		\ifx\IPmu\empty
			\varphi_{\FBstepsize*}^{\text{\sc fb}}
		\else
			\varphi_{\IPmu,\FBstepsize*}^{\text{\sc fb}}
		\fi
	}

	\newcommand{\LL}{\@ifstar\@@LL\@LL}
	\newcommand{\@LL}{\@ifnextchar_\@LL@sub\@LL@nosub}
	\newcommand{\@@LL}{\@ifnextchar_\@@LL@sub\@@LL@nosub}
	\newcommand{\@LL@nosub}{\mathcal L_{\nicefrac{1}{\FBstepsize}}}
	\newcommand{\@@LL@nosub}{\mathcal L_{\nicefrac{1}{\FBstepsize*}}}
	\def\@LL@sub_#1{\mathcal L_{#1}}
	\def\@@LL@sub_#1{\mathcal L_{#1}}

	\newcommand{\XX}{\R^n}
	\newcommand{\YY}{\R^m}

	\usepgfplotslibrary{groupplots}%
	\usetikzlibrary{matrix}%
	\definecolor{cb2blue}{RGB}{55,126,184}%
	\definecolor{cb2red}{RGB}{228,26,28}%
	\definecolor{cb2green}{RGB}{77,175,74}%
	\definecolor{cb2purple}{RGB}{152,78,163}%
	\definecolor{cb2orange}{RGB}{255,127,0}%
	\colorlet{lightblue}{cb2blue!40}%
	\colorlet{lightred}{cb2red!40}%
	\definecolor{mycoloreps1c5}{RGB}{252,174,145}%
	\definecolor{mycoloreps2}{RGB}{251,106,74}%
	\definecolor{mycoloreps2c5}{RGB}{222,45,38}%
	\definecolor{mycoloreps3}{RGB}{165,15,21}%
	\definecolor{mycolorregret1}{RGB}{186,228,188}%
	\definecolor{mycolorregret1c1}{RGB}{123,204,196}%
	\definecolor{mycolorregret1c25}{RGB}{67,162,202}%
	\definecolor{mycolorregret1c5}{RGB}{8,104,172}%
	\colorlet{mycolorregretAdapt}{cb2red}%

	\usetikzlibrary{shapes.geometric}
	\pgfdeclareplotmark{mystar1}{\node[draw=colorLimit1,star,star points=5,star point ratio=0.4,solid,fill=white,scale=.8] {};}
	\pgfdeclareplotmark{mystar2}{\node[draw=colorLimit2,star,star points=5,star point ratio=0.4,solid,fill=white,scale=.8] {};}
	\pgfdeclareplotmark{mystar3}{\node[draw=colorLimit3,star,star points=5,star point ratio=0.4,solid,fill=white,scale=.8] {};}

	\def\myampersand{&}

\makeatother

\newcommand{\TheTitle}{An interior proximal gradient method\\for nonconvex optimization}

\newcommand{\TheKeywords}{%
	Nonsmooth nonconvex optimization\AND%
	interior point methods\AND%
	proximal algorithms\AND%
	locally Lipschitz gradient%
}
\newcommand{\TheAMSsubj}{%
	\amsmscLink{49J52}\AND%
	\amsmscLink{65K05}\AND%
	\amsmscLink{90C30}%
}
\newcommand{\TheFunding}{%
	A. Themelis acknowledges the support of the Japan Society for the Promotion of Science (JSPS) KAKENHI grant JP21K17710.
}
\newcommand{\TheAddressADM}{%
	University of the Bun\-des\-wehr Munich\NL
	Department of Aerospace Engineering\NL
	Institute of Applied Mathematics and Scientific Computing\NL
	Werner-Heisenberg-Weg 39, 85577 Neubiberg, Germany%
}
\newcommand{\TheAddressAT}{%
	Kyushu University\NL
	Faculty of Information Science and Electrical Engineering (ISEE)\NL
	744 Motooka, Nishi-ku, 819-0395 Fukuoka, Japan%
}
\newcommand{\TheCodeZenodoDOI}{10.5281/zenodo.6890045}

	\author{%
		Alberto De Marchi\thanks{%
			\TheAddressADM\NL[.]%
			\emailLink{alberto.demarchi@unibw.de}\NL
			\orcidLink{0000-0002-3545-6898}%
		}%
		\and
		Andreas Themelis\thanks{%
			\TheAddressAT\NL[.]%
			\emailLink{andreas.themelis@ees.kyushu-u.ac.jp}\NL
			\orcidLink{0000-0002-6044-0169}.
			\TheFunding
		}%
	}
	\date{}
	\title{\TheTitle}

\let\includetikz\includegraphics
\graphicspath{{Pics/Tikz/}}

\begin{document}
	\maketitle
	\begin{abstract}
		We consider structured minimization problems subject to smooth inequality constraints and present a flexible algorithm that combines interior point (IP) and proximal gradient schemes.
		While traditional IP methods cannot cope with nonsmooth objective functions and proximal algorithms cannot handle complicated constraints, their combined usage is shown to successfully compensate the respective shortcomings.
		We provide a theoretical characterization of the algorithm and its asymptotic properties, deriving convergence results for fully nonconvex problems, thus bridging the gap with previous works that successfully addressed the convex case.
		Our interior proximal gradient algorithm benefits from warm starting, generates strictly feasible iterates with decreasing objective value, and returns after finitely many iterations a primal-dual pair approximately satisfying suitable optimality conditions.
		As a byproduct of our analysis of proximal gradient iterations we demonstrate that a slight refinement of traditional backtracking techniques waives the need for upper bounding the stepsize sequence, as required in existing results for the nonconvex setting.
	\end{abstract}
	\subclass{\TheAMSsubj}
	\tableofcontents

	\section{Introduction}
		We consider structured minimization problems
		\begin{mybox}
			\[
				\tag{P}\label{eq:P}
				\minimize_{x\in\R^n}\quad q(x) \coloneqq f(x)+g(x) \qquad
				\stt\quad c(x) \leq 0,
			\]
		\end{mybox}
		where \(\func{f}{\XX}{\R}\) and \(\func{c}{\XX}{\YY}\) are continuously differentiable and \(\func{g}{\XX}{\R\cup\set\infty}\) has easily computable proximal mapping.
		The structured objective \(q \coloneqq f+g\) is allowed to be nonconvex, as well as each component $f$ and $g$, and the constraint function \(c\) can be nonlinear.
		When the set induced by \(c(x) \leq 0\) is ``simple'', one may lift the inequality constraints to the objective of \eqref{eq:P}, enforcing them via an indicator function.
		But in many cases, projection onto the constraint set \(\set{x\in\XX}[c(x)\leq 0]\) can be expensive to compute, and even more so when coupled with the proximal mapping of \(g\), motivating us to seek a method able to handle inequalities explicitly.

		Starting from polynomial algorithms for linear programming \cite{khachiyan1979polynomial,karmarkar1984new}, interior point (IP) methods have shaken up the field of mathematical optimization and continue to spark renewed interest; see \cite{wright1997primal,forsgren2002interior,wright2005interior,gondzio2012interior} for a historical overview.
		It started by solving linear optimization problems with a nonlinear programming technique, based on the use of a barrier function \cite{frisch1955logarithmic} and sequential unconstrained minimization \cite{fiacco1968nonlinear}.
		The remarkable practical success was soon corroborated by deeper understanding of the major role played by the logarithmic barrier function \cite{gill1986projected,nesterov1994interior}, and similar methodologies were applied to solve quadratic and nonlinear optimization problems \cite{altman1999regularized,vanderbei1999interior,waechter2006implementation,curtis2012penalty,armand2017mixed}.
		However, the focus has almost exclusively been on smooth optimization and gradient-based or Newton-type methods.
		Some recent exceptions are the works on derivative-free \cite{brilli2022interior} and Riemannian \cite{lai2022riemannian} interior point methods for constrained optimization problems, as well as a closely related proximal gradient-based method \cite{chouzenoux2020proximal}.

		Recalling the basic idea of introducing a barrier function, the reader should observe that the IP rationale is independent of the smoothness of the functions defining the problem.
		Analogously to penalty and augmented Lagrangian methods \cite[\S 4.1]{birgin2014practical}, this feature contributes to the \emph{spirit of unification} that followed the interior point revolution \cite{forsgren2002interior}.
		But as far as we are aware, only a few articles consider IP approaches in the context of nonsmooth optimization problems such as \eqref{eq:P}.

		The combination of IP and splitting methods has been discussed by Valkonen \cite{valkonen2019interior} for a class of saddle point problems, associated with structured problems in the form $\min f + g \circ A$, where both $f$ and $g$ are possibly nonsmooth but convex, and $A$ is a bounded linear operator.
		More closely related to our approach, and associated with \eqref{eq:P}, is the proximal interior point algorithm (PIPA) presented in \cite{chouzenoux2020proximal}.
		Other works that depart from the classical Newton-type IP approach include \cite{lin2021admm}, which focuses on linear programs,
		and \cite{yang2023solving}, which addresses convex-constrained variational inequalities involving monotone operators.
		These works focus, however, on the convex setting and are not directly applicable if any of the problem data functions is nonconvex.
		Our work aims at filling this gap in the literature by developing and analyzing an interior point method for nonsmooth nonconvex problems.
		By extending the combination of splitting and IP methods to the fully nonconvex setting, we aim at bringing together and binding areas of optimization that seemed unrelated there.

		The constraint smoothening enabled by the adoption of suitably regular barriers in \eqref{eq:P} results in IP-type subproblems that seemingly retain a structure that proximal gradient iterations can address, namely the sum of a differentiable and a prox-friendly function.
		Seemingly, for both components are, in general, extended real-valued:
		the barrier term smoothens the (indicator of the) feasible set from the interior, thereby shrinking the domain of the differentiable term, as opposed to penalty (or augmented Lagrangian) schemes where the constraints are relaxed and the feasible set enlarged.
		Although sufficiently small stepsizes can be chosen to make gradient steps remain in the differentiable region, the composition with proximal operations precludes this possibility.
		Unless different techniques to deal with constraints are proposed, additional structural assumptions to prevent pathological instances are necessary.
		In the proximal interior point algorithm (PIPA) of \cite{chouzenoux2020proximal}, convexity is the key.

		Dropping these convexity assumptions, this work aims to be a first step toward wider applicability and more versatile modeling.
		In particular, we show that mere continuity of \(g\) \emph{relative to its domain} is sufficient, with no convexity restriction on any term of \eqref{eq:P}.
		This is achieved by leveraging an adaptive strategy that enables the use of proximal gradient both in absence of convexity and global Lipschitz differentiability requirements \cite{demarchi2022proximal,kanzow2022convergence}.
		With a detailed analysis around boundary points, where the barriers escape to infinity, \emph{local} properties are exploited to prove well definedness of the backtracking search.
		Then, we demonstrate that adaptive proximal gradient steps can generate (strictly) feasible iterates while guaranteeing a descent-type condition at the same time, eventually yielding an approximate KKT-optimal output.
		When specialized to the case \(c=0\) in \eqref{eq:P}, yet without \(g\) being necessarily continuous relative to its domain, it is shown that through a minor modification of the backtracking strategy no artificial bound on the stepsize sequence is necessary to recover standard convergence results for proximal gradient iterations, cf. \cref{thm:PG}.
		To the best of our knowledge, boundedness of the stepsize sequence is a standing assumption of any existing work dealing with the nonconvex case.

		We also point out the usage of non-Euclidean geometries induced by Bregman distances as another proximal gradient-based alternative to account for ambient constraints \cite{bolte2018first,lu2018relatively,latafat2022bregman,wang2023bregman}.
		Of this kind, Newton-type extensions also exist that can significantly speed up convergence and even attain superlinear rates, under assumptions at the limit point \cite{ahookhosh2021bregman,behmandpoor2022spiral}.
		All these methods are however subject to (and thus limited in applicability by) the identification of a distance-generating function enabling a so-called Lipschitz-like convexity condition, making induced proximal operations tractable, and whose domain agrees with the constraint set, which must thus be convex.
		Our focus is instead on addressing problem \eqref{eq:P} in the full generality of \cref{ass:basic}, stated next.

		\subsection{Problem setting and proposed methodology}
			We consider \eqref{eq:P} under the following standing assumptions.
			Technical definitions are given in \cref{sec:notation}.

			\begin{mybox*}
				\begin{ass}\label{ass:basic}%
					The following hold in problem \eqref{eq:P}:
					\begin{enumeratass}
					\item\label{ass:f}%
						\(\func{f}{\XX}{\R}\) has a locally Lipschitz-continuous gradient.
					\item\label{ass:g}%
						\(\func{g}{\XX}{\R\cup\set\infty}\) is proper, lsc, \(\gamma_g\)-prox-bounded, and continuous relative to \(\dom g\).
					\item\label{ass:c}%
						\(\func{c}{\XX}{\YY}\) has locally Lipschitz-continuous Jacobian.
					\item\label{ass:phi}%
						\(\inf\set{q(x)}[c(x) \leq 0 ] \in \R\).
					\item\label{ass:D}%
						The problem is strictly feasible: namely,
						\(\dom q\cap\set{x \in \XX}[c(x)<0]\neq\emptyset\).
					\end{enumeratass}
				\end{ass}
			\end{mybox*}
			From a computational point of view, it is assumed that one strictly feasible point can be retrieved explicitly, and that \(g\) has an easily computable proximal mapping.
			Continuity of \(g\) relative to its domain is meant in the sense that whenever \(\dom g\ni x^k\to x\) it holds that \(g(x^k)\to g(x)\).
			Few exceptions apart, such as functions involving 0-norms, most nonsmooth functions widely used in practice comply with this requirement.
			For instance, \(g\) can be the indicator of any nonempty and closed set, and thus enforce arbitrary closed constraints.

			The IP framework builds upon a barrier function $\func{b}{\R}{\R\cup\set\infty}$ to replace the \emph{inequality} constraints \cite{frisch1955logarithmic,fiacco1968nonlinear}.
			We will henceforth fix a nonnegative and smooth barrier function \(b\) that complies with the following requirements, assumed throughout.
			\begin{mybox*}
				\begin{ass}\label{ass:b}%
					The \emph{barrier function} $\func{b}{\R}{[0,\infty]}$ is such that
					\begin{enumeratass}
					\item
						$\dom b = (-\infty,0)$.
					\item
						$b$ is twice continuously differentiable with \(b'>0\) on its domain.
					\item
						$b(t) \to \infty$ as $t \to 0^-$.
					\end{enumeratass}
				\end{ass}
			\end{mybox*}
			Equality constraints should be considered carefully and treated e.g. via penalty \cite[\S 4.1.4]{curtis2012penalty} or augmented Lagrangian \cite{demarchi2023constrained} approaches.
			In the spirit of IP methods \cite{frisch1955logarithmic,fiacco1968nonlinear,bertsekas1999nonlinear,waechter2006implementation}, we consider a sequence of ``unconstrained'' barrier problems
			\begin{mybox}
				\begin{equation}
					\tag{P$\!_\mu$}\label{eq:Pmu}
					\minimize_{\x\in\R^n}\quad q_\mu(\x) \coloneqq f_\mu(\x)+g(\x),
				\end{equation}
			\end{mybox}
			whose differentiable cost function $\func{f_\mu}{\XX}{\R}$ includes the barrier terms weighted by a barrier parameter $\mu > 0$:%
			\begin{equation}\label{eq:fmu}
				f_\mu(\x) \coloneqq f(\x) + \mu \sum_{i=1}^m b( c_i(\x) ).
			\end{equation}
				The presence of the possibly nonsmooth term \(g\) prevents the employment of traditional IP methods which address the barrier subproblems by means of (smooth) Newton-type techniques.
				Instead, whenever \(g\) has an easily computable proximal mapping, instances of \eqref{eq:Pmu} are well suited for proximal gradient--based solvers.
				This is the rationale originally pursued in \cite{chouzenoux2020proximal} and that we here further extend beyond convexity assumptions.

			The procedure detailed in \cref{alg:IP} advances by minimizing the cost function at each iteration and updating the barrier parameter between iterations.
			At \cref{state:IP:x} a point \(x^{k+1}\) is retrieved by invoking the proximal gradient method \ipfb, outlined in \cref{alg:FB}, that provides a suitable numerical routine for addressing this task.
			Its definition requires some preliminary material and the introduction of some notation, and is therefore deferred to \cref{sec:FB}.
			The iterates $\seq{y^k}$ defined by \cref{state:IP:y} are solely involved in the termination criterion; as we will show, they relate to the Lagrange multipliers associated with the inequality constraints; cf. \cref{sec:preliminaries:stationarity}.

			\cref{alg:IP} provides a flexible template of an IP method for inequality constrained problems.
			It features warm-starting, inexact subsolves, and is subsolver-agnostic, meaning that one can run specialized routines for the problem at hand.
			In this work we focus on the proximal gradient-based \ipfb{} (\cref{alg:FB}), shown to be a suitable candidate for arbitrary formulations as \eqref{eq:P} whenever the proximal mapping of \(g\) is easily computable.

			\begin{algorithm}[t]
				\caption{\begin{tabular}[t]{@{}l@{}}
					Interior point method for \eqref{eq:P}\\
					using \ipfb{} (\cref{alg:FB}, page \pageref{alg:FB}) as inner subsolver%
				\end{tabular}}%
				\label{alg:IP}%
				\begin{algorithmic}[1]%
				\item[]\hspace*{-\leftmargin}%
				\begin{tabular}[t]{@{}l@{~~}l@{~~}l@{}}
					\algfont{Require} &
						\(x^0\) & strictly feasible starting point (\ie, \(x^0\in\dom g\) with \(c(x^0)<0\))
						\\
						& \(\epsilon_{\rm p},\epsilon_{\rm d}>0\) & primal-dual tolerances
				\\
					\algfont{Provide} &
						\(x^\star\)
						&
						\((\epsilon_{\rm p},\epsilon_{\rm d})\)-KKT optimal point for \eqref{eq:P} (cf. \cref{def:epsKKT})
				\\
					\algfont{Initialize} &
						\(\varepsilon_0,\mu_0>0\) \myampersand initial tolerance and barrier parameters
						\\
						& \(\theta_\varepsilon,\theta_\mu\in(0,1)\) \myampersand tolerance and barrier update coefficients
				\end{tabular}
				\hrule
				\vspace{5pt}%
				\linespread{1.1}\selectfont%
				\item[\algfont{repeat for} \(k=0,1,2\ldots\)]%
					\State\label{state:IP:x}%
						\(x^{k+1}=\ipfb(x^k,\mu_k,\varepsilon_k)\)~
					\Comment{%
						\(\varepsilon_k\)-stationary for \(q_{\mu_k}\)
						(see \cref{thm:FB:return})%
					}%
					\State\label{state:IP:y}%
						Set $y_i^{k+1} \gets \mu_k b'( c_i(x^{k+1}) )$ for all $i$
					\If{~
						\(\varepsilon_k \leq \epsilon_{\rm d}\)
						~\algfont{and}~
						\(\max_{i=1,\dots,m} \min\set*{ - c_i(x^{k+1}), y_i^{k+1} } \leq \epsilon_{\rm p}\)
					~}\label{state:IP:exit}
						\Statex*\algfont{return} \((x^\star,y^\star)\gets(x^{k+1},y^{k+1})\)
					\EndIf
					\State\label{state:IP:epsmu}%
						Select~
						$0 < \varepsilon_{k+1} \leq \max \set{\epsilon_{\rm d}, \theta_\varepsilon \varepsilon_k}$
						~and~
						$0 < \mu_{k+1} \leq \theta_\mu \mu_k$
				\end{algorithmic}
			\end{algorithm}

		\subsection{Contribution}
			We present an interior point proximal method (\cref{alg:IP}) for addressing in\-equal\-i\-ty-constrained structured minimization problems.
			Relying on suitable barrier functions and avoiding the need for slack variables to treat inequalities, our algorithm deviates from those based on penalty-type schemes \cite{sopasakis2020open,demarchi2023constrained}, and always generates feasible iterates while reducing the objective value.
			Convergence is guaranteed from arbitrary strictly feasible starting points (cf. \cref{thm:KKT,thm:IP:return}).
			To our knowledge, this work offers the first (feasible) IP method for addressing problem \eqref{eq:P} in the fully nonconvex setting.

			As a certified solver for the IP inner subproblems, we propose \ipfb, a proximal gradient method capable of handling barrier problems, whose well definedness is guaranteed through a suitable linesearch (cf. \cref{thm:FB:finite}).
			We establish convergence guarantees in the full generality of problems \eqref{eq:Pmu} (cf. \cref{thm:FB:asymp,thm:FB:return}), coping in particular with the lack of full domain of the smooth function therein.
			As a byproduct of our analysis, in \cref{thm:PG} we present the first convergence result of proximal gradient iterations with backtracking linesearch in a fully nonconvex regime that does not require any bound on the generated stepsize sequence.

		\subsection{Notation and known facts}\label{sec:notation}
			With $\N$, $\R$, \(\R_+\coloneqq[0,\infty)\) and $\Rinf \coloneqq \R \cup \set{\infty}$ we denote the natural, real, positive real, and extended-real numbers,
			respectively.
			Given $p \in \XX$ and a nonempty set $E \subset \XX$, \(\dist(p,E) \coloneqq \inf\set{\|x-p\|}[x\in E]\) denotes the distance of \(p\) from \(E\).
			The closed ball of radius \(r\) centered at \(p\) is denoted as \(\cball{p}{r}\coloneqq\set{x}[\|x-p\|\leq r]\).
			For a sequence \(\seq{x^k}\) and a set of indices \(K\subseteq\N\), \(x^k\to_Kx\) indicates that the subsequence \(\seq{x^k}[k\in K]\) converges to \(x\).

			Let $\func{F}{A}{\R^m}$ be a function defined on a set $A \subseteq \R^n$, and $\bar x \in A$.
			Following \cite[Def. 9.1]{rockafellar1998variational}, we say that $F$ is \DEF{locally Lipschitz} (or \DEF{strictly}) \DEF{continuous at $\bar x$} if $\bar x \in \interior A$ and the value
			\begin{equation}
				\lip F(\bar x) \coloneqq \limsup_{\substack{x,x^\prime \to \bar x\\x \ne x^\prime}} \frac{\|F(x) - F(x^\prime)\|}{\|x-x^\prime\|}
			\end{equation}
			is finite; here, \(\lip F(\bar x)\) denotes the Lipschitz constant of \(F\) at \(\bar x\).

			The notation \(\ffunc{T}{\R^n}{\R^n}\) indicates a point-to-set operator \(T\) that maps each \(x\in\R^n\) into a set \(T(x)\subseteq\R^n\).
			The \DEF{domain} of \(T\) is \(\dom T\coloneqq\set{x\in\R^n}[T(x)\neq\emptyset]\), and we say that \(T\) is \DEF{outer semicontinuous (osc)} if its \DEF{graph} \(\graph T\coloneqq\set{(x,y)}[y\in T(x)]\) is a closed subset of \(\R^n\times\R^n\).
			\(T\) is said to be \DEF{locally bounded} if for any bounded set \(E\subset\R^n\) it holds that \(\bigcup_{x\in E}T(x)\) is bounded.
			For a set-valued mapping, we use the \(\limsup\) notation to indicate the \DEF{outer limit} \cite[Def. 4.1]{rockafellar1998variational}, namely
			\[
				\bar y\in\limsup_{x\to\bar x}T(x)
			\quad\defeq[\Longleftrightarrow]\quad
				\exists\seq{x^k,y^k}\subseteq\graph T:
				\ (x^k,y^k)\to(\bar x,\bar y).
			\]
			In particular, \(T\) is osc if and only if \(T(\bar x)=\limsup_{x\to\bar x}T(x)\) for all \(\bar x\in\R^n\).

			The \DEF{effective domain} of an extended real-valued function $\func{h}{\R^n}{\Rinf}$ is denoted by $\dom h \coloneqq \set{x\in\R^n}[h(x) < \infty]$.
			We say that $h$ is \DEF{proper} if $\dom h \neq \emptyset$
			and \DEF{lower semicontinuous} (lsc) if $h(\bar x) \leq \liminf_{x\to\bar x} h(x)$ for all $\bar x \in \R^n$.
			For some constant $\tau\in\R$, $\lev_{\leq \tau} h\coloneqq\set{x\in\XX}[h(x)\leq \tau]$
			denotes the $\tau$-\DEF{sublevel set} associated with $h$.
			Following \cite[Def. 8.3]{rockafellar1998variational} and \cite[\S1.3]{mordukhovich2018variational}, we denote by $\ffunc{\hat{\partial} h}{\XX}{\XX}$ the \DEF{regular subdifferential} of $h$, where
			\begin{equation}
				\bar v \in \hat{\partial} h(\bar x)
			\quad\defeq[\Leftrightarrow]\quad
				\liminf_{\substack{x\to\bar x\\x\neq\bar x}} \frac{h(x) - h(\bar x) - \innprod{\bar v}{x-\bar x}}{\|x-\bar x\|} \geq 0 .
			\end{equation}
			The (\DEF{limiting}) \DEF{subdifferential} of $h$ is $\ffunc{\partial h}{\XX}{\XX}$, where $\bar v \in \partial h(\bar x)$ if and only if \(\bar x\in\dom h\) and there exist sequences $\seq{x^k}$ and $\seq{v^k}$ such that $(x^k,v^k,h(x^k))\to(\bar x,\bar v, h(\bar x))$ and $v^k \in \hat\partial h(x^k)$ for all \(k\).
				By considering a constant sequence \(x^k\equiv\bar x\), the inclusion \(\hat\partial h(\bar x)\subseteq\partial h(\bar x)\) readily follows.
						The subdifferential of $h$ at $\bar x$ satisfies $\partial(h+h_0)(\bar x) = \partial h(\bar x) + \nabla h_0(\bar x)$ for any $\func{h_0}{\XX}{\Rinf}$ continuously differentiable around $\bar x$ \cite[Ex. 8.8]{rockafellar1998variational}.

			The \DEF{proximal mapping} of \(h\) with stepsize \(\gamma>0\) is the set-valued operator \(\ffunc{\prox_{\gamma h}}{\R^n}{\R^n}\) defined as%
			\begin{equation}\label{eq:prox}
				\prox_{\gamma h}(x)
			\coloneqq
				\argmin_{z\in\R^n} \set{ h(z) + \tfrac{1}{2\gamma}\|z-x\|^2},
			\end{equation}
			and we say that $h$ is \DEF{prox-bounded} if it is proper and $h + \frac{1}{2\gamma}\|{}\cdot{}\|^2$ is bounded below on $\XX$ for some $\gamma > 0$.
			The supremum of all such $\gamma$ is the \DEF{threshold $\gamma_h$ of prox-boundedness} for $h$.
			In particular, if $h$ is bounded below by an affine function, then $\gamma_h = \infty$.
			When $h$ is lsc, for any $\gamma \in (0,\gamma_h)$ and \(x\in\R^n\) it holds that \cite[Thm 1.25]{rockafellar1998variational}
			\begin{equation}\label{eq:proxOSC}
				\emptyset
			\neq
				\limsup_{(x',\gamma')\to(x,\gamma)}\prox_{\gamma'h}(x')
			\subseteq
				\prox_{\gamma h}(x).
			\end{equation}

	\section{Stationarity and optimality concepts}\label{sec:preliminaries:stationarity}
		Iterative minimization methods typically approach local solutions only asymptotically, while in finitely many iterations can only yield points that satisfy some relaxed, or approximate, optimality conditions.
		In the case of the minimization of a proper function \(\func{h}{\R^n}{\Rinf}\), the inclusion \(0\in\partial h(x^\star)\) (in fact, \(0\in\hat\partial h(x^\star)\)) is necessary for local minimality of \(x^\star\) for \(h\) \cite[Thm 10.1]{rockafellar1998variational}.
		An approximate counterpart can be formulated by bounding the distance of the zero vector from the subdifferential.
		The following definition introduces a terminology tailored for inner problem instances \eqref{eq:Pmu}.

		\begin{mybox}
		\begin{defin}[\(\varepsilon\)-stationarity for \eqref{eq:Pmu}]%
			Relative to \eqref{eq:Pmu}, a point \(x^\star\) is \DEF{\(\varepsilon\)-stationary} for some \(\varepsilon\geq0\) if
			\[
				\dist(0,\partial  q_\mu(x^\star))\leq\varepsilon.
			\]
			When \(\varepsilon=0\), \ie, when \(0\in \partial  q_\mu(x^\star)\), \(x^\star\) is said to be \DEF{stationary}.\footnotemark
		\end{defin}
		\end{mybox}

		\footnotetext{%
			The equivalence of \(\dist(0,\partial  q_\mu(x^\star))=0\) and \(0\in \partial  q_\mu(x^\star)\) follows from closedness of \(\partial  q_\mu(x^\star)\), see \cite[Thm. 8.6]{rockafellar1998variational}.
		}%
		Considering the minimization problem defining the proximal mapping as in \eqref{eq:prox}, the necessary stationarity condition reads
		\begin{equation}\label{eq:proxsubgrad}
			\tfrac{x-\bar x}{\gamma}
		\in
			\hat\partial h(\bar x)
		\subseteq
			\partial h(\bar x)
		\qquad
			\forall \bar x\in\prox_{\gamma h}(x).
		\end{equation}
		Notice that whenever \(x^\star\) is an (approximate) stationary point for \eqref{eq:Pmu}, it necessarily belongs to the domain of \(q_\mu\), for otherwise \(\partial  q_\mu(x^\star)\) would be empty.
		In particular, \(c(x^\star)<0\), a stronger condition than that prescribed by the constraint in the original problem \eqref{eq:P}.
		To emphasize the difference, we will talk in terms of feasibility and \emph{strict} feasibility, as defined next.

		\begin{mybox}
			\begin{defin}[Strict feasibility]\label{def:feasible}%
				Relative to problem \eqref{eq:P}, a point $x^\star \in \dom q$ is called \DEF{feasible} if $c(x^\star) \leq 0$, and \DEF{strictly feasible} if \(c(x^\star)<0\).
			\end{defin}
		\end{mybox}

		The given notion of (strict) feasibility imposes the inclusion \(x^\star\in\dom q\) so as to also account for implicit constraints encoded in the cost function.
		Problem \eqref{eq:P} can equivalently be expressed as the ``unconstrained'' minimization of the extended real-valued function  		%
		\begin{equation}\label{eq:q0}
			q_0\coloneqq q+\indicator_{\R_-^m}\circ c,
		\end{equation}
		where for a set \(E\subseteq\YY\) we denote by \(\func{\indicator_E}{\YY}{\Rinf}\) the \DEF{indicator function} of \(E\), defined as \(\indicator_E(x)=0\) if \(x\in E\) and \(\infty\) otherwise.
		In these terms, feasibility of \(x^\star\) can be expressed as the inclusion \(x^\star\in\dom q_0\), whereas strict feasibility as the inclusion \(x^\star\in\dom q_\mu\) for some (in fact, any) \(\mu>0\).
		The notion of feasibility is therefore independent of how the problem is formulated, whereas the set of strictly feasible points depends on the specific representation of $g$ and $c$.

		Similarly, in addressing problem \eqref{eq:P} one could in principle seek for (approximate) stationary points of \(q_0\).
		In practice, however, complications may arise in resolving the nonsmooth subdifferential chain rule involved in the evaluation of \(\partial  q_0\).
		For this reason, following the nonlinear programming approach we will consider KKT-type optimality conditions when dealing with \eqref{eq:P}.
		These constitute a relaxed stationarity condition, and are in fact equivalent under suitable constraint and epigraphical qualifications.

		\begin{mybox}
			\begin{defin}[KKT optimality for \eqref{eq:P}]\label{def:KKT}
				Relative to \eqref{eq:P}, a point $x^\star\in\R^n$ is \DEF{KKT-optimal} if it is feasible and there exists $y^\star \in \R_+^m$ such that
				\begin{subequations}\label{subeq:KKT}%
					\begin{gather}\label{eq:KKT:x}
						-\nabla c(x^\star)^\top y^\star
					\in
						\partial q(x^\star)
					\shortintertext{and}
						y_i^\star c_i(x^\star) = 0
						\quad \forall i=1,\dots,m.
						\label{eq:KKT:y}
					\end{gather}
				\end{subequations}
			\end{defin}
		\end{mybox}

		Mirroring the concept of \(\varepsilon\)-stationarity for ``unconstrained'' minimization problems such as \eqref{eq:Pmu}, the next definition gives a characterization of approximate KKT optimality for problems subject to (explicit) constraints.
		This notion allows us to qualify the output of \cref{alg:IP} in relation to \eqref{eq:P}; similarly, approximate stationarity will serve as the counterpart for the ``unconstrained'' inner subproblems \eqref{eq:Pmu}.

		\begin{mybox}
			\begin{defin}[\((\epsilon_{\rm p},\epsilon_{\rm d})\)-KKT optimality for \eqref{eq:P}]\label{def:epsKKT}%
				Relative to \eqref{eq:P}, a point $x^\star\in\R^n$ is said to be \DEF{\((\epsilon_{\rm p},\epsilon_{\rm d})\)-KKT optimal} for some \(\epsilon_{\rm p},\epsilon_{\rm d}\geq0\) if it is feasible and there exists $y^\star \in \R_+^m$ such that
				\begin{subequations}\label{subeq:epsKKT}%
					\begin{gather}\label{eq:epsKKT:x}
						\dist\bigl(-\nabla c(x^\star)^\top y^\star,\,\partial q(x^\star)\bigr)
					\leq
						\epsilon_{\rm d}
					\shortintertext{and}
						\min\set{-c_i(x^\star), y_i^\star}
					\leq
						\epsilon_{\rm p}
						\quad \forall i=1,\dots,m.
						\label{eq:epsKKT:y}
					\end{gather}
				\end{subequations}
			\end{defin}
		\end{mybox}

		Notice that, together with feasibility of $x^\star$ and nonnegativity of $y^\star$, condition \eqref{eq:epsKKT:y} imposes a constraint of approximate complementarity.
		In general, it is not weaker nor stronger than the more classical condition $|y_i^\star c_i(x^\star)| \leq \epsilon_{\rm p}$, which could be considered as well.

		Similarly to what remarked for approximate stationarity, \((\epsilon_{\rm p},\epsilon_{\rm d})\)-KKT optimality naturally reduces to KKT optimality when \(\epsilon_{\rm p}=\epsilon_{\rm d}=0\).
		There is, however, a substantial difference in the behavior of approximate stationary and approximate KKT optimal points when the tolerances approach zero in the limit.
		Suppose that \(\seq{\x^k}\) is an \(\varepsilon_k\)-stationary point for \eqref{eq:Pmu}, with \(\varepsilon_k\searrow0\) and \(\x^k\to\x^\star\).
		Under \cref{ass:basic}, we may immediately deduce that \(\x^\star\) is stationary.\footnote{%
			In absence of continuity of \(g\) on its domain, the claim still holds true provided that \(\x^k\) converges \(q_\mu\)-attentively, namely in such a way that \(q_\mu(\x^k)\to q_\mu(\x^\star)\).%
		}
		On the contrary, having \(x^k\to x^\star\) with \((\epsilon_{{\rm p},k},\epsilon_{{\rm d},k})\)-KKT optimal for \eqref{eq:P} and \(\epsilon_{{\rm p},k},\epsilon_{{\rm d},k}\searrow 0\) does not guarantee KKT optimality of the limit \(x^\star\).
		This issue raises the need of explicitly defining an asymptotic version of approximate KKT optimality, on the vein of \cite[Def. 3.1]{birgin2014practical} and \cite[Def. 2.4]{demarchi2023constrained}.

		\begin{mybox}
			\begin{defin}[A-KKT optimality]\label{def:AKKT}%
				Relative to \eqref{eq:P}, a point $x^\star\in\R^n$ is said to be \emph{asymptotically KKT (A-KKT) optimal} if it is feasible and there exist \(\seq{y^k}\subset\R_+^m\) and a feasible sequence $\seq{x^k}\to x^\star$ such that%
				\begin{subequations}\label{subeq:AKKT}
					\begin{gather}\label{eq:AKKT:x}%
						\dist\bigl(-\nabla c(x^k)^\top y^k,\,\partial q(x^k)\bigr)\to0
					\shortintertext{and}
					\label{eq:AKKT:y}%
						y_i^kc_i(x^\star)=0
						\quad\forall i=1,\dots,m.
					\end{gather}
				\end{subequations}
			\end{defin}
		\end{mybox}

		Having \(y_i^kc_i(x^\star)=0\) in condition \eqref{eq:AKKT:y} causes no loss of generality over \(y_i^kc_i(x^\star)\to0\), a seemingly more natural asymptotic counterpart of \eqref{eq:KKT:y}.
		This equivalence will be useful in the sequel, and is formally stated in the following lemma for future reference.

		\begin{mybox}
			\begin{lem}\label{thm:yk}%
				Suppose that \cref{ass:basic} holds, and let a feasible sequence \(\seq{x^k}\subset\R^n\) converging to a feasible point \(x^\star\) and a sequence \(\seq{\tilde y^k}\subset\R_+^m\) be such that
				\begin{subequations}\label{subeq:AKKT'}
					\begin{gather}\label{eq:AKKT':x}%
						\dist\bigl(-\nabla c(x^k)^\top\tilde y^k,\,\partial q(x^k)\bigr)\to0
					\shortintertext{and}
					\label{eq:AKKT':y}
						\tilde y_i^kc_i(x^\star)\to0
						\quad\forall i=1,\dots,m.
					\end{gather}
				\end{subequations}
				Then, \(x^\star\) is A-KKT optimal.
			\end{lem}
		\end{mybox}
		\begin{proof}
			For all \(k\in\N\) and \(i=1,\dots,m\), define \(y_i^k=\tilde y_i^k\) if \(c_i(x^\star)=0\) and \(y_i^k=0\) otherwise.
			Then, observing that \(\|\tilde y_i^k-y_i^k\|\to0\), it is immediate to verify that \(\seq{x^k}\) and \(\seq{y^k}\) comply with \cref{def:AKKT}.
		\end{proof}

		It is also worth remarking that usual notions of A-KKT optimality do not require feasibility of the points \(x^k\); nevertheless, in our setting where these points are retrieved through inner IP procedures, feasibility (in fact, \emph{strict}) comes at no cost since it is always inherently satisfied.

		While KKT clearly implies A-KKT, the discrepancy between the two notions is again to be found in unmet qualifications, in absence of which local minimizers may fail to be KKT optimal, even for convex problems; A-KKT optimality, on the contrary, is necessary.
		In referring the reader to the well documented \cite[\S3]{birgin2014practical} for examples and a thorough discussion, we point out that the feature of A-KKT optimality allowing it to encompass any local solution lies in the possible unboundedness of the sequence \(\seq{y^k}\) in \cref{def:AKKT}, in absence of which the notion reduces to the nonasymptotic KKT counterpart.

		\begin{rem}\label{thm:ykbounded}%
			If the sequence \(\seq{y^k}\) in \cref{def:AKKT} has a cluster point \(y^\star\), as is the case when it is bounded, then the point \(x^\star\) therein is KKT optimal, not only asymptotically.
			This simply follows from the continuity of \(q\) on its domain, implying that \(\limsup_{k\to\infty}\partial q(x^k)\subseteq\partial q(x^\star)\), and hence that
			\begin{align*}
				0
			={} &
				\lim_{k\to\infty}\dist\bigl(-\nabla c(x^k)^\top y^k,\,\partial q(x^k)\bigr)
			\\
			\geq{} &
				\limsup_{k\to\infty}\dist\bigl(-\nabla c(x^k)^\top y^k,\,\partial q(x^\star)\bigr)
			\\
			={} &
				\dist\bigl(-\nabla c(x^\star)^\top y^\star,\,\partial q(x^\star)\bigr)
			\end{align*}
			by continuity of \(\nabla c\) and of the distance function.
		\end{rem}

\def\iter{j}{}%
	\section{A barrier-friendly proximal gradient method}\label{sec:FB}
		In this section we elaborate upon \cref{state:IP:x} of \cref{alg:IP}, that aims at solving the barrier problem \eqref{eq:Pmu} via proximal gradient iterations.
		Specifically, we will show that at every (outer) iteration \(k\), the call to \ipfb{} yields a point \(x^{k+1}\) which is \(\varepsilon_k\)-stationary for problem (\hyperref[eq:Pmu]{\(P_{\mu_k}\)}) and such that \(q_{\mu_k}(x^{k+1})\leq q_{\mu_k}(x^k)\), as commented at \cref{state:IP:x}.
		\ipfb, outlined in \cref{alg:FB}, is adapted from \cite[Alg. 3]{demarchi2022proximal} so as to cope with the lack of the full domain of the locally smooth function \(f_\mu\).
		In fact, improving upon \cite{demarchi2022proximal,kanzow2022convergence,demarchi2023proximal} we here remove boundedness impositions on the stepsize sequence.
		This flexibility is captured, at the beginning of every iteration \(j\), by initializing the stepsize as \(\gamma_j=r\gamma_{j-1}\) (as opposed to \(\gamma_j=\gamma_{j-1}\), or selecting \(\gamma_j\) from a fixed bounded interval),
		where the factor \(r\geq1\) quantifies the stepsize enlargement.
		Large values of \(r\) aim at expediting convergence in terms of number of iterations by testing large stepsizes first, at the expense of potentially more backtrackings and, consequently, gradient evaluations per iteration.
		Small values instead result in fewer backtrackings at the expense of more conservative stepsize choices.
		By compensating for the possibly overly cautious estimate obtained by previous reductions, this stepsize redemption has been denominated \emph{``regret''} in the FOM toolbox \cite{beck2019fom}, a terminology that we also adopt in this work.
		Although the tuning of \(r\) may be problem dependent,
		recent results for the convex case provide insights on parameter-free and problem-independent choices;
		we refer to the commentary after \cref{thm:PG} for the details.

		Relative to \eqref{eq:Pmu}, we consider the proximal gradient operator with stepsize $\gamma\in(0,\gamma_g)$ defined by
		\begin{equation}
			\def\iter{}%
			\T(\x) \coloneqq \FB{\x}
			\def\iter{j}%
		\end{equation}
		which is compact valued, and relative to
		\[
			\dom \T_{\mu,\gamma} = \dom f_\mu = \set{\x \in \XX}[c(\x) < 0]
		\]
		it is outer semicontinuous (osc) and locally bounded.\footnote{\label{footnote:Tosc}%
			Local boundedness relative to \(\dom f_\mu\) indicates that for every compact set \(Z\subset\dom f_\mu\) the set \(\bigcup_{\x\in Z}\T_{\mu,\gamma}(\x)\) is bounded.
			Moreover, for any \(z\in\dom f_\mu\) and \(\gamma\in(0,\gamma_g)\) it follows from \eqref{eq:proxOSC} that
			\(
				\emptyset
			\neq
				\limsup_{(\x',\gamma')\to(\x,\gamma)}\T_{\mu,\gamma'}(\x')
			\subseteq
				\T_{\mu,\gamma}(\x)
			\).%
		}
		Notice that, in general, the range of \(\T_{\mu,\gamma}\) need not be contained in its domain; as such, fixed-point iterations of \(\T_{\mu,\gamma}\) may be ill defined.

		\begin{algorithm}[t]
			\caption{\begin{tabular}[t]{@{}l@{}}
				\protect\ipfb\((\x,\mu,\varepsilon)\)\\
				\def\textbf{}%
				\textbf Forward \textbf Backward solver for \textbf Inner \textbf Problem \eqref{eq:Pmu}%
			\end{tabular}}%
			\label{alg:FB}%
			\begin{algorithmic}[1]%
			\newcommand{\inlinegap}{\hspace*{0.15cm}}%
			\item[]\hspace*{-\leftmargin}%
				\begin{tabular}[t]{@{}l@{~~}l@{~~}l@{}}
					\algfont{Require}
						& \(\x\)  & strictly feasible starting point (\ie, \(\x\in\dom g\) with \(c(\x)<0\))\\
						& \(\mu>0\)         & barrier coefficient\\
						& \(\varepsilon>0\) & termination tolerance
				\\
					\algfont{Provide}
						& \(\x^\ast\) & (strictly feasible) \(\varepsilon\)-stationary point for \eqref{eq:Pmu}
				\\
					\algfont{Initialize} &
						\(\gamma_0\in (0,\gamma_g)\) & initial stepsize
						\\
						& \(\alpha,\beta\in(0,1)\) & stepsize backtracking parameters
						\\
						& \(r\geq1\) & stepsize regret factor
				\end{tabular}
			\linespread{1.1}\selectfont%
			\item[Set \(\x^0\gets\x\) and \algfont{repeat for \(j=0,1,\dots\)}]%
				\hrule
				\vspace{5pt}%

			\State \label{state:FB:init}\label{state:FB:x+}%
				\algfont{if}\inlinegap
					\(j\geq1\)
				\inlinegap\algfont{then}\inlinegap
					\(\gamma_j\gets r\gamma_{j-1}\) and \(\x^j\gets\bar\x^{j-1}\);
				\inlinegap\algfont{end if}
				\Comment{(\(\gamma_j\gets\min\set{r\gamma_{j-1},\gamma_g-\delta}\) for some \(\delta>0\) if \(\gamma_g\neq\infty\))}

			\While{ \algfont{true} }\label{state:FB:while}%

				\State \label{state:FB:barx}%
					Compute $\bar\x^j \in \T(\x^j)$%
					\vspace*{4pt}%

				\State\label{state:FB:gammaLS}%
					\algfont{if}\inlinegap
						\(\left\{\begin{minipage}{6.5cm}%
							\begin{conditions}
							\item \label{cond:FB:domfmu}%
								\(c(\bar\x^j)<0\)
							\item \label{cond:FB:geq}%
								\(
									q_\mu(\bar\x^j)
									\leq
									q_\mu(\x^j) - \tfrac{1-\alpha}{2\gamma_j}\|\bar\x^j - \x^j\|^2
								\)
							\item \label{cond:FB:Lip}%
								\(\|\nabla f_\mu(\bar\x^j)-\nabla f_\mu(\x^j)\|\leq\frac{\alpha}{\gamma_j}\|\bar\x^j-\x^j\|\)
							\end{conditions}%
						\end{minipage}~\right\}\)
					\inlinegap\algfont{then}\inlinegap
						{\algfont{break}};
					\inlinegap\algfont{else}\inlinegap
						$\gamma_j \gets \beta \gamma_j$;
					\inlinegap\algfont{end if}%

			\EndWhile

				\State\label{state:FB:exit}%
					\algfont{if}\inlinegap
						\( \| \frac{1}{\gamma_j} (\x^j - \bar\x^j) - \nabla f_\mu(\x^j) + \nabla f_\mu(\bar\x^j) \| \leq \varepsilon\)~
					\inlinegap\algfont{then}\inlinegap
						{\algfont{return}} \(\x^\ast \gets \bar\x^j\)
					\inlinegap\algfont{end if}%
			\end{algorithmic}
		\end{algorithm}

		Beyond the introduction of the regret factor \(r\), the results and proofs stated in the following closely pattern those presented in \cite{demarchi2022proximal}, where proximal gradient with an adaptively tuned stepsize is shown to work under a mere local Lipschitz differentiability assumption of the smooth term.
		Although \cref{alg:FB} is effectively a classical adaptive proximal gradient method, the challenge here is twofold.
		First, the range of the proximal gradient operator may fail to be contained in its domain, which precludes the possibility of a na\"ive fixed-point approach.
		Second, the adaptive strategy considered in \cite{demarchi2022proximal} revolves around the fact that in any bounded set a finite modulus of Lipschitz continuity of the gradient of the smooth function exists; this property dramatically fails for \(f_\mu\) in the IP setting here investigated, as its gradient explodes whenever approaching the boundary of the constraint set \(\set{\x \in \XX}[c(\x)\leq0]\).
		While these issues have been examined and well resolved in \cite{chouzenoux2020proximal} for the convex case, no successful attempt appears to have been accomplished in the nonconvex setting.

		The key difference with traditional proximal gradient settings is that here, under \cref{ass:basic}, the function \(f_\mu\) as defined in \eqref{eq:fmu} has (locally) Lipschitz-continuous gradient \emph{on its domain}, as opposed to on the entire space.
		This means that for every \emph{convex and compact} set \(\Omega\subset\dom f_\mu\) there exists \(L_{f_\mu,\Omega}\geq0\) such that
		\begin{equation}\label{eq:fmuC1+}
			\begin{cases}
				\|\nabla f_\mu(\x')-\nabla f_\mu(\x)\|\leq L_{f_\mu,\Omega}\|\x'-\x\|
				\\[3pt]
				f_\mu(\x')\leq f_\mu(\x)+\innprod{\nabla f_\mu(\x)}{\x'-\x}+\tfrac{L_{f_\mu,\Omega}}{2}\|\x'-\x\|^2
			\end{cases}
			\quad
			\forall \x,\x'\in\Omega,
		\end{equation}
		see \cite[Thm. 9.2]{rockafellar1998variational} and \cite[Prop. A.24]{bertsekas1999nonlinear}.
		In fact, as detailed in the former reference, one can take
		\(
			L_{f_\mu,\Omega}
		=
			\sup_{\x\in\Omega}\lip\nabla f_\mu(\x)
		\)
		in this case.
		Nevertheless, an elementary compactness argument shows that a finite \(L_{f_\mu,\Omega}\) exists for any compact \emph{but not necessarily convex} \(\Omega\subset\dom f_\mu\).
		This observation suggests that, inasmuch as the iterates are confined sufficiently far away from the troublesome boundary of \(\set{z}[c(z)\leq0]\), issues originating from the lack of full domain of \(f_\mu\) can be circumvented.
		A simple proof for the validity of \eqref{eq:fmuC1+} for any compact \(\Omega\subset\dom f_\mu\) is detailed for completeness.
		Note that the interpretation of \(L_{f_\mu,\Omega}\) as a Lipschitz constant is ill posed when the set \(\Omega\) is not convex, and the supremum formula only furnishes a lower bound to \(L_{f_\mu,\Omega}\) in this case.

		\begin{mybox}
		\begin{lem}\label{thm:fmuC1+}%
			Let \(\mu>0\) be fixed.
			For any compact set \(\Omega\subset\dom f_\mu=\set{\x}[c(\x)<0]\) there exists a constant \(L_{f_\mu,\Omega}\geq0\) satisfying \eqref{eq:fmuC1+}.
		\end{lem}
		\end{mybox}
		\begin{proof}
			Contrary to the claim, suppose that for any \(j\in\N\) there exist \(\x_j,\x'_j\in\Omega\) violating either one of the two conditions in \eqref{eq:fmuC1+} with \(L_{f_\mu,\Omega}=j\) therein.
			By compactness of \(\Omega\), there exists an infinite index set \(J\subseteq\N\) together with \(\x,\x'\in\Omega\) such that \(\x_j\to\x\) and \(\x'_j\to\x'\) as \(J\ni j\to\infty\).
			Since \(\x,\x'\in\Omega\subset\dom f_\mu\) and \(\Omega\) is compact, necessarily \(\x=\x'\) (for otherwise finiteness of either \(f_\mu(\x)\), \(f_\mu(\x')\), \(\nabla f_\mu(\x)\), or \(\nabla f_\mu(\x')\) would be violated).
			As a consequence, up to discarding early terms if necessary openness of \(\dom f_\mu\) entails the existence of \(\delta>0\) such that \(\x_j,\x'_j\in\cball\x\delta\subset\dom f_\mu\) holds for all \(j\in J\).
			This is a contradiction, since for any \(j\geq L_{f_\mu,\cball{\x}{\delta}}\) both conditions hold, where the existence of \(L_{f_\mu,\cball{\x}{\delta}}\geq0\) is guaranteed by compactness \emph{and convexity} of \(\cball{\x}{\delta}\subset\dom f_\mu\).
		\end{proof}

		\subsection{Algorithm outline}
			Although retaining the core features of the adaptive proximal gradient method \cite[Alg. 3]{demarchi2022proximal}, see Corollary 4.7 therein, \ipfb{} includes checks in order to generate iterates that are strictly feasible for $c(\x) \leq 0$ and exhibit a sufficient decrease on the cost function.
			These conditions are enforced at \cref{state:FB:gammaLS}; notice that \cref{cond:FB:domfmu} is implied by \cref{cond:FB:geq}, and could thus be safely removed without affecting the algorithm.
			We however prefer to explicitly include the former as well both for clarity and algorithmic convenience: assessing \cref{cond:FB:geq} requires evaluating \(c(\bar{\x}^j)\) in the first place and, if \cref{cond:FB:domfmu} is found to fail, the whole if statement can already be resolved to be false without further unnecessary function evaluations.
			Notice further that since $\bar{\x}^j = \x^{j+1}$, $\nabla f_\mu(\bar{\x}^j)$ evaluated within the \(j\)-th iteration can be stored and used in the next one to save computations.

			Finite termination of the linesearch occurring at \cref{state:FB:gammaLS} hinges on the strict feasibility of the previous iterate, which is why the condition must be satisfied in the first place by the initial point $\x^0$ fed in input to \cref{alg:FB}.
			When called within the IP routine of \cref{alg:IP} at \cref{state:IP:x}, this condition is always inherently satisfied, since the initial point \(x^k\) prescribed therein is the ouput of a previous call to \ipfb, and is thus strictly feasible by construction.
			By estimating the local Lipschitz constant of $\nabla f_\mu$ and monitoring the cost function $q_\mu$, the algorithm is shown to generate iterates $\seq{\bar{\x}^j}[j\in\N]$ that remain bounded away from the barrier at $\set{\x\in\XX}[c(\x) = 0]$.
			As mentioned in the foreword to \cref{thm:fmuC1+}, this is the key feature to circumvent the lack of full domain of \(f_\mu\).

			As will be shown in \cref{thm:FB:return}, the termination criterion at \cref{state:FB:exit} is satisfied in finitely many iterations and entails \(\varepsilon\)-stationarity of the output \(\bar\x^j\) for \(q_\mu\).
			The condition is clearly satisfied if \(\bar\x^j=\x^j\), in which case \(\bar\x^j\) is stationary, not only approximately so.
			For this reason, without loss of generality we may avoid trivialities by assuming throughout that $\bar{\x}^j \neq \x^j$ holds for every \(j\).

		\subsection{Well definedness}
			We start by observing that each problem instance \eqref{eq:Pmu} is well posed, and also list some important structural properties as placeholders for future reference.
			The proof of the assertions is a trivial consequence of \cref{ass:basic,ass:b}.

			\begin{mybox}
				\begin{lem}\label{thm:FB:properties}%
					For any $\mu > 0$, the following hold:
					\begin{enumerate}
					\item \label{thm:FB:properties:phi}%
						$\func{q_\mu}{\XX}{\Rinf}$ is proper, lsc, with \(\dom q_\mu=\dom g\cap\set{\x}[c(\x)<0]\) and \(\inf q_\mu\in\R\).
					\item \label{thm:FB:properties:psi}%
						$\func{f_\mu}{\XX}{\Rinf}$ has locally Lipschitz gradient on \(\dom f_\mu=\set{\x}[c(\x)<0]\).
					\end{enumerate}
				\end{lem}
			\end{mybox}

			We proceed to show that \ipfb{} is well defined, namely that each iteration successfully terminates without getting stuck in infinite loops at \cref{state:FB:while}.
			Our argument is based on the fact that the proximal mapping converges to the identity as the stepsize tends to zero, a claim that is formalized in the following auxiliary result.

			\begin{mybox}
			\begin{lem}\label{thm:xbartox}%
				Let \(\func{h}{\R^n}{\Rinf}\) be lsc, and \(\seq{\x^\ell}[\ell\in\N]\subset\R^n\) be a sequence converging to a point \(\x\in\dom h\).
				Let \(\bar\x^\ell\in\prox_{\gamma_\ell h}(\x^\ell)\) with \(\gamma_\ell\searrow0\).
				Then, \(\bar\x^\ell\to\x\).
			\end{lem}
			\end{mybox}
			\begin{proof}
				We start by observing that the existence of \(\bar\x^\ell\) guarantees prox-bound\-ed\-ness (hence properness) of \(h\).
				For every \(\ell\), the optimality of \(\bar\x^\ell\) in the proximal minimization subproblem reads
				\(
					h(\bar\x^\ell)+\tfrac{1}{2\gamma_\ell}\|\bar\x^\ell-\x^\ell\|^2
				\leq
					h(\x)+\tfrac{1}{2\gamma_\ell}\|\x-\x^\ell\|^2
				\).
				By invoking the triangle and Young's inequalities, this implies that
				\begin{align*}
					\| \bar\x^\ell - \x \|^2
				\leq{} &
					\| \bar\x^\ell - \x^\ell \|^2 + 2 \| \bar\x^\ell - \x^\ell \| \| \x - \x^\ell \| + \| \x - \x^\ell \|^2
				\\
				\leq{} &
					2 \| \bar\x^\ell - \x^\ell \|^2 + 2 \| \x - \x^\ell \|^2
				\\
				\leq{} &
					4 \gamma_\ell h(\x) + 2 \|\x-\x^\ell\|^2 - 4 \gamma_\ell h(\bar\x^\ell) + 2 \| \x - \x^\ell \|^2
				\\
				={} &
					4 \left[ \gamma_\ell h(\x) - \gamma_\ell h(\bar\x^\ell) + \| \x - \x^\ell \|^2 \right]
					.
				\end{align*}
				By rearranging, we obtain
				\begin{equation}\label{eq:proxoptim}
					\gamma_\ell h(\bar\x^\ell)+\tfrac{1}{4}\|\bar\x^\ell-\x\|^2
				\leq
					\gamma_\ell h(\x)
				+
					\|\x-\x^\ell\|^2.
				\end{equation}
				The right-hand side vanishes as \(\ell\to\infty\); since \(h\) is proper and lsc, it suffices to show that \(\seq{\bar\x^\ell}[\ell\in\N]\) remains bounded, as this would imply that each term on the left-hand side too vanishes as \(\ell\to\infty\).
				Contrary to the claim, up to extracting, suppose that \(\|\bar\x^\ell\|\to\infty\).
				Then, dividing both sides of \eqref{eq:proxoptim} by \(\|\bar\x^\ell\|^2\) yields
				\[
					\liminf_{\ell\to\infty}\gamma_\ell\tfrac{h(\bar\x^\ell)}{\|\bar\x^\ell\|^2}
				\leq
					-\tfrac14,
				\quad\text{hence}\quad
					\liminf_{\ell\to\infty}\tfrac{h(\bar\x^\ell)}{\|\bar\x^\ell\|^2}
				=
					-\infty.
				\]
				By virtue of \cite[Ex. 1.24]{rockafellar1998variational}, this contradicts prox-boundedness of \(h\).
			\end{proof}

			\begin{mybox}
				\begin{lem}[Well definedness]\label{thm:FB:finite}%
					Consider \eqref{eq:Pmu} and the iterates generated by \cref{alg:FB}.
					The following hold:
					\begin{enumerate}
					\item\label{thm:FB:finite:LS}%
						At every iteration, the number of backtrackings at \cref{state:FB:gammaLS} is finite.
					\item\label{thm:FB:finite:descent}\label{thm:FB:x+=barx}%
						At the \(j\)-th iteration (\(j\geq1\)),
						one has
						\(
							\x^j=\bar\x^{j-1}
						\)
						and
						\begin{equation}\label{eq:FB:SD}
							q_\mu({\x}^j)
						=
							q_\mu(\bar{\x}^{j-1})
						\leq
							q_\mu(\x^{j-1})
							-
							\tfrac{1-\alpha}{2\gamma_{j-1}}
							\|\bar \x^{j-1}-\x^{j-1}\|^2.
						\end{equation}
					\item\label{thm:FB:finite:sublevel}%
						Every iterate \(\bar \x^j\) remains within \(\lev_{\leq q_\mu^0} q_\mu\), where \(q_\mu^0 \coloneqq q_\mu(\x^0) < \infty\).
					\end{enumerate}
				\end{lem}
			\end{mybox}
			\begin{proof}
				Let us index by \({j,\ell}\) the variables defined at the \(\ell\)-th attempt within the $j$-th iteration.
				\begin{proofitemize}
				\item{thm:FB:finite:LS}
					Let us show that from some strictly feasible $\x^{j-1}$, $j\geq1$, the iteration terminates (in finite time) yielding a strictly feasible $\x^j$.
					Terminating an iteration requires to satisfy the conditions at \cref{state:FB:gammaLS}.
					To arrive to a contradiction, suppose that this never happens, hence that \(\gamma_{j,\ell} = \beta^{\ell}r\gamma_{j-1} \searrow0\) as \(\ell\to\infty\).
					By openness of \(\dom f_\mu\ni \x^{j-1}\), there exists \(\delta_j>0\) such that \(\Omega_j\coloneqq\cball{\x^{j-1}}{\delta_j}\subset\dom f_\mu\).
					Since \(\x^{j-1} - \gamma_{j,\ell} \nabla f_{\mu}(\x^{j-1}) \to \x^{j-1} \in \dom g\) as \(\gamma_{j,\ell} \searrow 0\), \cref{thm:xbartox} applies and yields the existence of \(\ell_j\geq0\) such that \(\bar\x^{j,\ell}\in\Omega_j\) for all \(\ell\geq \ell_j\).
					On the other hand, by convexity and compactness of \(\Omega_j\subset\dom f_\mu\), for any given $\alpha \in (0,1)$ there also exists \(\ell_j'\geq0\) such that \(\nicefrac{\alpha}{\gamma_{j,\ell}}\geq L_{f_\mu,\Omega_j}\) for all \(\ell\geq \ell_j'\).
					From \cref{thm:fmuC1+} we then conclude that for any \(\ell\geq\max\set*{\ell_j,\ell_j'}\) both conditions at \cref{state:FB:gammaLS} are satisfied.
					In particular, for \(\ell\geq\max\set*{\ell_j,\ell_j'}\) we have%
					\[
						f_\mu(\bar{\x}^{j,\ell})
					\leq
						f_\mu(\x^{j,\ell})
						+ \innprod*{\nabla f_\mu(\x^{j,\ell})}{\bar{\x}^{j,\ell} - \x^{j,\ell}}
						+ \tfrac{\alpha}{2\gamma_{j,\ell}}\|\bar{\x}^{j,\ell} - \x^{j,\ell}\|^2 .
					\]
					Meanwhile, the minimizing property of $\bar{\x}^{j,\ell}$ at \cref{state:FB:barx} implies
					\[
						g(\bar{\x}^{j,\ell})
						+
						\innprod*{\nabla f_\mu(\x^{j,\ell})}{\bar{\x}^{j,\ell} - \x^{j,\ell}}
						+
						\tfrac{1}{2\gamma_{j,\ell}}\|\bar{\x}^{j,\ell} - \x^{j,\ell}\|^2
					\leq
						g(\x^{j,\ell}).
					\]
					Combining these inequalities, \cref{cond:FB:geq} is eventually satisfied, whence the contradiction.

				\item{thm:FB:finite:descent}
					The assertion follows from the failure of the condition at \cref{state:FB:gammaLS} and the fact that the value of \(\x^j\) is not updated after its definition at \cref{state:FB:init}.%

				\item{thm:FB:finite:sublevel}
					Follows from assertion \ref{thm:FB:finite:descent}, with $q_\mu(\x^0) < \infty$ since $\x^0$ is strictly feasible.
				\qedhere
				\end{proofitemize}
			\end{proof}

		\subsection{Convergence analysis}
			The remainder of the section is devoted to showing that for every strictly feasible initial point \(\x\) and \(\mu,\varepsilon>0\) \ipfb\((\x,\mu,\varepsilon)\) returns an \(\varepsilon\)-stationary point \(\x^\star\) for \(q_\mu\) satisfying \(q_\mu(\x^\star)\leq q_\mu(\x)\).
			To this end, we provide an asymptotic analysis where we show that with \(\varepsilon=0\) the algorithm runs indefinitely and produces iterates satisfying
			\(
				\liminf_{j\to\infty}\| \frac{1}{\gamma_j} (\x^j - \bar{\x}^j) - \nabla f_\mu(\x^j) + \nabla f_\mu(\bar{\x}^j) \|
			=
				0
			\),
			see \cref{thm:FB:asymp:res}.
			The claimed successful finite termination can then be deduced, as ultimately formalized in \cref{thm:FB:return}.
			The entire proof of \cref{thm:FB:asymp} is carried out without assuming continuity of \(g\) on its domain as required in \cref{ass:g}.
			In allowing the stepsize regret parameter \(r\) to be strictly greater than 1, and without imposing any upper bound on the stepsizes \(\gamma_j\) (other than staying bounded away from \(\gamma_g\), should this threshold be finite), this theorem constitutes an important refinement of \cite[Cor. 4.7]{demarchi2022proximal} and other related works on proximal gradient algorithms such as \cite{salzo2017variable,kanzow2022convergence,demarchi2023proximal} which rely on boundedness of \(\seq{\gamma_j}[j\in\N]\).

			\begin{mybox}
				\begin{thm}[Asymptotic analysis of \ipfb]\label{thm:FB:asymp}%
					The iterates generated by \cref{alg:FB} with termination tolerance \(\varepsilon=0\) satisfy the following:%
					\begin{enumerate}[itemsep=0pt,topsep=2pt]
					\item\label{thm:FB:asymp:cost}%
						$\seq{ q_\mu(\x^j) }[j\in\N]$ converges to a finite value $q_\mu^\star \geq \inf q_\mu$ from above.

					\item\label{thm:FB:asymp:summable}%
						$\sum_{j\in\N} \frac{1}{\gamma_j}\|\bar \x^j-\x^j\|^2 < \infty$.

					\item\label{thm:FB:asymp:feas}%
						\(\sup_{j\in\N}\max\set*{c_i(\bar \x^j),c_i(\x^j)}<0\), for every \(i=1,\dots,m\).

					\item\label{thm:FB:asymp:gammaconstant}%
						Consider the following assertions:
						\begin{enumerate}[label={(\oldstylenums{\arabic*})},ref={(\oldstylenums{\arabic*})},topsep=0pt]
						\item\label{thm:FB:gammaconstant:lb}%
							$q_\mu$ is level bounded;

						\item\label{thm:FB:gammaconstant:barxk}%
							$\seq{ \bar{\x}^j }[j\in\N]$ is bounded;

						\item\label{thm:FB:gammaconstant:xk}%
							$\seq{ \x^j }[j\in\N]$ is bounded;

						\item\label{thm:FB:gammaconstant:gamma}%
							$\seq{ \gamma_j }[j\in\N]$ is bounded away from zero, \ie,
							there exists $\gamma_{\rm min} > 0$ such that $\gamma_j \geq\gamma_{\rm min}$ for every $j$.%
						\end{enumerate}
						One has~
						\ref{thm:FB:gammaconstant:lb}
						~$\Rightarrow$~
						\ref{thm:FB:gammaconstant:barxk}
						~$\Leftrightarrow$~
						\ref{thm:FB:gammaconstant:xk}
						~$\Rightarrow$~
						\ref{thm:FB:gammaconstant:gamma}.

					\item\label{thm:FB:asymp:gamma}%
						$\sum_{j\in\N} \gamma_j = \infty$.

					\item\label{thm:FB:asymp:res}%
						\(\displaystyle
							\liminf_{j\to\infty} \tfrac{1}{\gamma_j}\| \bar{\x}^j - \x^j \|
						=
							\liminf_{j\to\infty}\| \tfrac{1}{\gamma_j} (\x^j - \bar{\x}^j) - \nabla f_\mu(\x^j) + \nabla f_\mu(\bar{\x}^j) \|
						=
							0
						\).
					\item\label{thm:FB:subseq}
						If the iterates remain bounded, then the set \(\omega\) of accumulation points of $\seq{ \bar{\x}^j }[j\in\N]$ is made of stationary points for \(q_\mu\), and \(q_\mu\) is constantly equal to \(q_\mu^\star\) as in assertion \ref{thm:FB:asymp:cost} on \(\omega\).
					\end{enumerate}
					All these claims hold without \(g\) being necessarily continuous relative to its domain.
				\end{thm}
			\end{mybox}
			\begin{proof}%
				We begin by observing that (the proofs of) all the claims of \cref{thm:FB:properties,thm:xbartox,thm:FB:finite} that we shall refer to hereafter are indipendent of whether \(g\) is continuous on its domain or not.
				\begin{proofitemize}
				\item{thm:FB:asymp:cost}
					Follows from \cref{thm:FB:finite:descent,thm:FB:properties:phi}.

				\item{thm:FB:asymp:summable}
					Follows from a telescoping argument on \eqref{eq:FB:SD}, having
					\begin{equation}\label{eq:FB:summable}
						(1-\alpha)
						\sum_{j\in\N}{
							\tfrac{1}{2\gamma_j}
							\|\bar \x^j-\x^j\|^2
						}
					\leq
						q_\mu(\x^0)
						-
						\inf q_\mu
					<
						\infty.
					\end{equation}

				\item{thm:FB:asymp:feas}
					Let \(i\in\set{1,\dots,m}\) be fixed.
					For every \(j\in\N\) we have
					\[
						\inf\set{q(\x)}[c(\x) \leq 0 ]
					+
						\mu b(c_i(\bar \x^j))
					\leq
						q(\bar \x^j)
					+
						\mu b(c_i(\bar \x^j))
					\leq
						q_\mu(\bar \x^j)
					\leq
						q_\mu(\x^0),
					\]
					where the infimum attains a finite value by \cref{ass:phi}, since \(b\geq0\),
					the second inequality too uses nonnegativity of \(b\), and the last one follows from \cref{thm:FB:finite:sublevel}.
					Therefore, the sequence \(\seq{b(c_i(\bar \x^j))}[j\in\N]\) remains bounded, which implies that \(\seq{c_i(\bar \x^j)}[j\in\N]\) is bounded away from 0.
					In turn, since \(\x^j=\bar\x^{j-1}\) by \cref{thm:FB:x+=barx}, so is \(\seq{c_i(\x^j)}[j\in\N]\).

				\item{thm:FB:asymp:gammaconstant}
					The first implication follows from \cref{thm:FB:finite:sublevel}, and the second one from \cref{thm:FB:x+=barx}.
					Suppose now that \(\seq{\x^j}[j\in\N]\) is bounded, and thus that so is \(\seq{\bar \x^j}[j\in\N]\).
					From assertion \ref{thm:FB:asymp:feas} we then infer the existence of a compact set \(\Omega\subset\dom f_\mu\) that contains both sequences.
					As argued in the proof of \cref{thm:FB:finite:LS}, any value \(\gamma_j\leq\nicefrac{\alpha}{L_{f_\mu,\Omega}}\) will pass all conditions at \cref{state:FB:gammaLS} and will thus not be subject to any backtracking.

				\item{thm:FB:asymp:gamma}
					By iteratively applying the triangle inequality (recall that \(\x^j=\bar\x^{j-1}\), cf. \cref{thm:FB:x+=barx}), we obtain%
					\begin{align*}
						\|\x^j-\x^0\|
					\leq{} &
						\sum_{\ell=0}^{j-1}\|\bar \x^\ell-\x^\ell\|
					=
						\sum_{\ell=0}^{j-1}\gamma_\ell^{-\nicefrac12}\|\bar \x^\ell-\x^\ell\|\gamma_\ell^{\nicefrac12}
					\\
					\leq{} &
						\sqrt{
							\sum_{\ell=0}^{j-1}
							\gamma_\ell^{-1}\|\bar \x^\ell-\x^\ell\|^2
						}
						\sqrt{
							\sum_{\ell=0}^{j-1}\gamma_\ell
						}
					\overrel*[\leq]{\eqref{eq:FB:summable}}{} 
						\sqrt{
							2 \frac{q_\mu(\x^0) - \inf q_\mu}{1-\alpha}
							\vphantom{\sum_{\ell=0}^{j-1}\gamma_\ell}
						}
						\sqrt{
							\sum_{\ell=0}^{j-1}\gamma_\ell
						}.
					\end{align*}
					Contrary to the claim, if \(\sum_{j\in\N}\gamma_j<\infty\) holds, then \(\seq{\x^j}[j\in\N]\) is bounded.
					From assertion \ref{thm:FB:asymp:gammaconstant} we then infer that \(\gamma_j\) is bounded away from zero, thus contradicting the finiteness of \(\sum_{j\in\N}\gamma_j\).

				\item{thm:FB:asymp:res}
					That \(\liminf_{j\to\infty} \tfrac{1}{\gamma_j}\| \bar{\x}^j - \x^j \|=0\) follows from assertions \ref{thm:FB:asymp:summable} and \ref{thm:FB:asymp:gamma}.
					In turn, the other limit follows from the fact that
					\(
						\|\nabla f_\mu(\x^j)-\nabla f_\mu(\bar{\x}^j)\|
					\leq
						\frac{\alpha}{\gamma_j}
						\| \bar{\x}^j - \x^j \|
					\),
					enforced by \cref{cond:FB:Lip}.

				\item{thm:FB:subseq}
					\def\FBstepsize{\gamma}%
					It follows from assertions \ref{thm:FB:asymp:feas} and \ref{thm:FB:asymp:gammaconstant} that the iterates \(\x^j\) and \(\bar\x^j\) are contained in a compact set \(\Omega\subset\dom f_\mu\), and that \(\gamma_j\geq\gamma_{\rm min}>0\) holds for all \(j\).
					Let \(\x^\star\in\omega\) be fixed and let an infinite set of indices $J\subseteq\N$  be such that $\bar{\x}^j \to_J \x^\star$.
					Observe that optimality of \(\bar\x^j\) in the minimization problem defining \(\T_{\mu,\gamma_{j}}(\x^j)\) implies
					\[
						g(\bar\x^j)
						+
						\tfrac{1}{2\gamma_{j}}\|\bar\x^j-\x^j+\gamma_{j}\nabla f_\mu(\x^j)\|^2
					\leq
						g(\x^\star)
						+
						\tfrac{1}{2\gamma_{j}}\|\x^\star-\x^j+\gamma_{j}\nabla f_\mu(\x^j)\|^2,
					\]
					which after expanding the squares and using the fact that \(\gamma_j\geq\gamma_{\rm min}>0\) gives
					\[
						g(\bar\x^j)
					\leq
						g(\x^\star)
						+
						\tfrac{1}{2\gamma_{\rm min}}
						\|
							\overbracket[0.5pt]{\x^\star-\bar\x^j}^{\to_J0}
						\|^2
						+
						\innprod*{
							\overbracket[0.5pt]{\nabla f_\mu(\x^j)}^{\text{\clap{bounded}}}
						}{
							\overbracket[0.5pt]{\x^\star-\bar\x^j}^{\to_J0}
						}
						-
						\tfrac{1}{2\gamma_j}\|\bar\x^j-\x^j\|^2.
					\]
					Therefore, \(\limsup _{J\ni j\to\infty}g(\bar\x^j)\leq g(\x^\star)\).
					Because of lsc, necessarily \(g(\bar\x^j)\to_J g(\x^\star)\), which together with continuity of \(f_\mu\) on \(\Omega\) leads to \(q_\mu(\bar\x^j)\to_J q_\mu(\x^\star)\).
					From the definition of \(q_\mu^\star\) in assertion \ref{thm:FB:asymp:cost} it then follows that \(q_\mu(\x^\star)=q_\mu^\star\), and the arbitrarity of $\x^\star \in \omega$ yields that \(q_\mu\equiv q_\mu^\star\) on \(\omega\).

					To prove stationarity, we consider two cases.
					If, up to extracting, \(\gamma_j\to_J\gamma<\gamma_g\leq\infty\), then the vanishing of \(\frac{1}{\gamma_j}\|\x^j-\bar \x^j\|^2\) implies that
					\[
						\x^\star
					=
						\lim_{J\ni j\to\infty}\bar\x^j
					\in
						\limsup_{J\ni j\to\infty}\T_{\mu,\gamma_j}(\x^j)
					\subseteq
						\T(\x^\star)
					\defeq
						\FB{\x^\star}
					\]
					with the last inclusion owing to outer semicontinuity of \(\T\) on \(\Omega\) (cf. \cref{footnote:Tosc}).
					The inclusion \(\x^\star\in\FB{\x^\star}\) together with \eqref{eq:proxsubgrad} yields the claimed stationarity \(0\in\hat\partial q_\mu(\x^\star)\subseteq\partial q_\mu(\x^\star)\).
					If, instead, \(\gamma_j\to_J\infty\), then since \(\x^j,\bar\x^j\) range in a bounded set,
					\(
						\|\nabla f_\mu(\x^j) - \nabla f_\mu(\bar{\x}^j)\|
					\leq
						\frac{\alpha}{\gamma_j}
						\| \bar{\x}^j - \x^j \|
					\to_J
						0
					\),
					where the first inequality is enforced at \cref{cond:FB:Lip}.
					It then follows that
					\(
						v^j
					\coloneqq
						\tfrac{1}{\gamma_j} (\x^j - \bar{\x}^j) - \nabla f_\mu(\x^j) + \nabla f_\mu(\bar{\x}^j)
					\to_J
						0
					\).
					Noticing that \(v^j\in\nabla f(\bar\x^j)+\hat\partial g(\bar\x^j)=\hat\partial q_\mu(\bar\x^j)\), cf. \eqref{eq:proxsubgrad}, and recalling that \(q_\mu(\x^j)\to_Jq_\mu(\x^\star)\) as shown above, we conclude that \(0\in\partial q_\mu(\x^\star)\).
				\qedhere
				\end{proofitemize}
			\end{proof}

			We can now easily infer finite termination of \ipfb{} for any \(\varepsilon>0\), confirming that the output of \ipfb{} is a feasible input for the outer IP framework of \cref{alg:IP}, as commented in \cref{state:IP:x} therein.

			\begin{mybox}
				\begin{cor}[\ipfb{} as inner solver for \cref{alg:IP}]\label{thm:FB:return}%
					For any strictly feasible starting point \(\x\) and \(\mu,\varepsilon>0\), in finitely many steps \ipfb\((\x,\mu,\varepsilon)\) returns an \(\varepsilon\)-stationary point \(\x^\star\) for \eqref{eq:Pmu} satisfying \(q_\mu(\x^\star)\leq q_\mu(\x)\).
				\end{cor}
			\end{mybox}
			\begin{proof}
				That the algorithm terminates in finitely many iterates, say \(j\) many, follows from \cref{thm:FB:asymp:res}.
				Since \(\bar\x^j\in\T(\x^j)=\FB{\x^j}\), it follows from \eqref{eq:proxsubgrad} that the output \(\x^\star=\bar\x^j\) satisfies
				\[
					\tfrac{1}{\gamma_j} (\x^j - \bar{\x}^j) - \nabla f_\mu(\x^j) + \nabla f_\mu(\bar{\x}^j)
				\in
					\hat\partial g(\bar{\x}^j) + \nabla f_\mu(\bar{\x}^j)
				=
					\hat\partial q_\mu(\bar{\x}^j)
				\subseteq
					\partial q_\mu(\bar{\x}^j)
				.
				\]
				The magnitude of such subgradient is no more than \(\varepsilon\) as enforced by the termination criterion, implying that \(\bar\x^j\) is \(\varepsilon\)-stationary for \(q_\mu\).
				Finally, that \(q_\mu(\bar\x^j)\leq q_\mu(\x)\) follows from \cref{thm:FB:finite:sublevel}.
			\end{proof}

			Incidentally, when specialized to the case \(c=0\), \cref{thm:FB:asymp} offers insights on plain proximal gradient (PG) iterations that, to the best of our knowledge, are novel.
			Specifically, it shows that enforcing a Lipschitz-like condition in addition to the standard quadratic upper bound allows one to waive any artificial cap on the stepsize sequence, which is a standing assumption in related literature.
			The chosen terminology \emph{``unconstrained stepsizes''} emphasizes this distinction.

			\begin{mybox}
				\begin{thm}[Convergence of PG with unconstrained stepsizes]\label{thm:PG}%
					Let \(\varphi\coloneqq f+g\) for a differentiable function \(\func{f}{\R^n}{\R}\) with locally Lipschitz-continuous gradient and a proper, lsc, and \(\gamma_g\)-prox bounded function \(\func{g}{\R^n}{\Rinf}\).
					Starting from \(\x^0\in\R^n\) and \(\gamma_0\in(0,\gamma_g)\), and given some \(\alpha,\beta\in(0,1)\) and \(r\geq1\), consider the following scheme:
					\begin{center}
						\fbox{%
							\begin{minipage}{0.8\linewidth}
								\normalfont\small
								\algfont{for}~~ \(j=1,2,\dots\) ~~\algfont{do}
								\begin{algsteps}
								\item
									\algfont{while ~~true~~ do}

									\item
										\qquad
										\(\x^j\in\prox_{\gamma_jg}(\x^{j-1}-\gamma_j\nabla f(\x^{j-1}))\)

									\item\label{state:PG:if}%
										\qquad
										\algfont{if}~~
										\(
											\varphi(\x^j)
										\leq
											\varphi(\x^{j-1})-\frac{1-\alpha}{2\gamma_j}\|\x^j-\x^{j-1}\|^2
										\)
										~~\algfont{and}~~
										\(
											\|\nabla f(\x^j)-\nabla f(\x^{j-1})\|
										\leq
											\frac{\alpha}{\gamma_j}\|\x^j-\x^{j-1}\|
										\)
										~~\algfont{then}

										\item
											\qquad
											\qquad
											\algfont{break}

									\item
										\qquad
										\(\gamma_j\gets\beta\gamma_j\)

								\item \label{state:PG:gamma_init}%
									\(\gamma_{j+1}=r\gamma_j\)
								\hfill
									{\color{gray}\(\triangleright\)\footnotesize\sf~(or \(\gamma_{j+1}=\min\set{r\gamma_j,\gamma_g-\delta}\) for some \(\delta>0\) in case \(\gamma_g\neq\infty\))}
								\end{algsteps}
							\end{minipage}
						}%
					\end{center}
					Then, \(\sum_{j\in\N}\gamma_j=\infty\) and \(\liminf_{j\to\infty}\bigl\|\frac{1}{\gamma_j}(\x^{j-1}-\x^j)-\bigl(\nabla f(\x^{j-1})-\nabla f(\x^{j-1})\bigr)\bigr\|=0\).
					If \(\seq{\x^j}[j\in\N]\) is bounded (\eg when \(\varphi\) is level bounded), then its cluster set \(\omega\) is made of stationary points for \(\varphi\), \(\varphi|_\omega\equiv\lim_{j\to\infty}\varphi(\x^j)\), and \(\inf_{j\in\N}\gamma_j>0\).
				\end{thm}
			\end{mybox}
			\begin{proof}
				We shall see this as a special case of \ipfb{} with \(c=0\) and \(\mu=0\), resulting in \(\dom f_\mu=\dom f=\R^n\) and thus with \cref{cond:FB:domfmu} vacuously satisfied at any backtracking test.
				If \(\x^0\notin\dom g\), then in the first iteration the first condition at \cref{state:PG:if} is also vacuously satisfied for (any candidate) iterate \(\x^1\).
				On the other hand, the second condition is satisfied for \(\gamma_j\) small enough, because of local Lipschitz continuity of \(\nabla f\) (and the fact that all the iterates \(\x^1\) tested in the backtracking remain in a bounded set).
				Then, for any \(j\geq1\) (regardless of whether \(\x^0\in\dom g\) or not) it holds that \(\x^j\in\dom g\), it being the output of a proximal mapping of \(g\).
				From iteration \(j=1\) on, then, we may invoke the proof of \cref{thm:FB:asymp}.
			\end{proof}

			Some comments are in order.
			The Lipschitz-like condition
			\(
				\|\nabla f(\x^j)-\nabla f(\x^{j-1})\|
			\leq
				\frac{\alpha}{\gamma_j}\|\x^j-\x^{j-1}\|
			\)
			at \cref{state:PG:if} in the PG scheme synopsized in \cref{thm:PG}, this being the refinement that allows for unbounded stepsizes, comes at a price, for every failed assessment incurs a wasted evaluation of \(\nabla f(\x^j)\).

			We also remark that the first condition
			\[
				\varphi(\x^j)
			\leq
				\varphi(\x^{j-1})-\frac{1-\alpha}{2\gamma_j}\|\x^j-\x^{j-1}\|^2
			\]
			is implied by the usual local quadratic upper bound
			\[
				f(\x^j)
			\leq
				f(\x^{j-1})+\innprod{\nabla f(\x^{j-1})}{\x^j-\x^{j-1}}+\frac{\alpha}{2\gamma_j}\|\x^j-\x^{j-1}\|^2,
			\]
			cf. the proof of \cref{thm:FB:finite:descent}.
			The validity of \cref{thm:PG} is thus unaffected if within the backtracking the latter inequality is adopted instead, which has the advantage of saving evaluations of \(g\) at the expense of a slight additional conservatism.

			Notice that the regret factor \(r\), that is, the ratio between the initial stepsize at any iteration and the accepted value at the previous one, is chosen constant for notational convenience and simplicity of exposition, but any sequence \(\seq{r_j}[j\in\N]\subset[1,\infty)\) would be an equally valid option.
			In other words, the stepsize initialization at \cref{state:PG:gamma_init} can be replaced by any \(\gamma_{j+1}\geq\gamma_j\) (as long as this choice is bounded away from \(\gamma_g\), should this threshold be finite).
			Nevertheless, a small parameter in the range \(r\in(1,2]\) is found to work particularly well in practice, an observation that recent results in the convex setting, advocating an adaptive \(r_j=\sqrt{1+\nicefrac{\gamma_{j-1}}{\gamma_{j-2}}}\), may shed some light upon; see \cite{malitsky2020adaptive} for the pioneering analysis in the smooth case and the follow-up proximal extensions \cite{latafat2023adaptive,malitsky2023adaptive,latafat2023convergence}, in particular the discussion surrounding \cite[Thm. 1]{malitsky2023adaptive}.
			This parallel is further emphasized by the Lipschitz-like condition
			\(
				\gamma_j\frac{\|\nabla f(\x^j)-\nabla f(\x^{j-1})\|}{\|\x^j-\x^{j-1}\|}
			\leq
				\alpha
			\),
			though the stepsize index is shifted in the cited references which allows one to waive any backtrack altogether in the convex case.
			In the analysis of \cref{thm:FB:asymp} and its special case \cref{thm:PG}, this Lipschitz-like condition is the key for lifting boundedness requirements on the stepsize sequence.

\def\iter{k}%
	\section{The outer interior point framework}\label{sec:IP}
		In the nonsmooth setting associated to \eqref{eq:P}, a proximal gradient algorithm such as \ipfb{} can be adopted for computing an approximate solution of subproblems in the form of \eqref{eq:Pmu}, as shown in \cref{sec:FB}.
		The choice of the first parameter (\ie, the initial point for the inner problem) in the call to \ipfb{} at \cref{state:IP:x} is dictated by the following rationale.
		Practical performances of both inner and outer procedure may benefit from warm-starting.
		The similarity between inner problem instances in subsequent iterations, namely instances of \eqref{eq:Pmu} solely differing by a slight variation of the parameter \(\mu\), suggests that the (approximate) solution \(x^k\) of the previous inner problem is an educated choice as initial iterate for the starting point of the current one.
		Furthermore, being the output of a call to \ipfb, \(x^k\) is guaranteed to be strictly feasible (for \(k=0\) this is true by initialization), and its employment as starting point for \ipfb{} is thus also theoretically supported.

		We proceed with a characterization of the iterates generated by \cref{alg:IP}, in terms of objective value, feasibility and stationarity.

		\begin{mybox}
			\begin{lem}[Algorithmic behavior]\label{thm:barriermethodprop}%
				Consider a sequence $\seq{x^k,y^k}$ generated by \cref{alg:IP}.
				For every \(k\geq0\), the following hold:
				\begin{enumerate}[itemsep=0pt]
				\item\label{thm:decrease}%
					\(q(x^{k+1}) \leq q_{\mu_k}(x^{k+1}) \leq q_{\mu_k}(x^k) \leq q_{\mu_{k-1}}(x^k)\).
				\item\label{thm:xkfeas}%
					\(x^{k+1}\) is \(\varepsilon_k\)-stationary for \(q_{\mu_k}\), and is in particular strictly feasible: \(x^{k+1}\in\dom q\) and \(c(x^{k+1})<0\).
				\item\label{thm:ykfeas}%
					\(y^{k+1}\geq0\).
				\item\label{thm:xkopt}%
					\(\dist\bigl(-\nabla c(x^{k+1})^\top y^{k+1},\partial q(x^{k+1})\bigr)\leq\varepsilon_k\).
				\end{enumerate}
			\end{lem}
		\end{mybox}
		\begin{proof}
			We remind that \(x^{k+1}\) is the output of \ipfb\((x^k,\mu_k,\varepsilon_k)\), cf. \cref{state:IP:x}.
			\begin{proofitemize}
			\item{thm:decrease}
				The second inequality follows from \cref{thm:FB:return}, and the other two from the fact that \(b\geq0\) and \(0\leq\mu_k\leq\mu_{k-1}\).

			\item{thm:xkfeas}
				Follows from \cref{thm:FB:return}.

			\item{thm:ykfeas}
				Follows from the fact that \(b'\geq0\) and \(\mu_k\geq0\).

			\item{thm:xkopt}
				$\varepsilon_k$-stationarity of \(x^{k+1}\) for \(q_{\mu_k}\) reads
				\(
					\dist(0,\partial q_{\mu_k}(x^{k+1}))\leq\varepsilon_k
				\).
				The claim then follows by observing that
				\begin{align*}
					\partial q_{\mu_k}(x^{k+1})
				={} &
					\partial q(x^{k+1})+\mu_k\sum_{i=1}^mb'(c_i(x^{k+1}))\nabla c_i(x^{k+1})
				\\
				={} &
					\partial q(x^{k+1})+\nabla c(x^{k+1})^\top y^{k+1},
				\end{align*}
				where the last identity uses the definition of \(y^{k+1}\) at \cref{state:IP:y}.
			\qedhere
			\end{proofitemize}
		\end{proof}

		We next turn our attention to finite termination and output qualification for \cref{alg:IP}.
		Similarly to the analysis carried out for the inner \ipfb{} in the previous section, we will obtain the results as a simple consequence of a more general asymptotic analysis in which the tolerances are driven to zero.

		\begin{mybox}
			\begin{thm}[Asymptotic analysis of \cref{alg:IP}]\label{thm:KKT}%
				Consider a sequence $\seq{x^k, y^k}$ of iterates generated by \cref{alg:IP}.
				Then,
				\begin{enumerate}[resume]
				\item\label{thm:KKT:bounded}%
					If the problem is coercive, in the sense that \(q_0\) as in \eqref{eq:q0} is level bounded, then \(\seq{x^k}\) is bounded.
				\item\label{thm:KKT:feas}%
					Any limit point of \(\seq{x^k}\) is feasible.
				\item\label{thm:KKT:ep}%
					If either \(\epsilon_{\rm p}=0\) or \(\epsilon_{\rm d}=0\), then \(\lim_{k\to\infty}\min\set{-c(x^k),y^k}=0\).
				\end{enumerate}
				If \(\epsilon_{\rm d}=0\), so that the algorithm runs indefinitely with \(\varepsilon_k,\mu_k\to0\), the following also hold for a subsequence $\seq{x_k}[k\in K]$ converging to a point \(x^\star\):
				\begin{enumerate}[resume]%
				\item\label{thm:KKT:AKKT}%
					$x^\star$ is a (feasible) A-KKT-optimal point for \eqref{eq:P}.
				\item\label{thm:KKT:KKT}%
					If $\seq{y^k}[k\in K]$ remains bounded, then \(x^\star\) is a KKT-optimal point for \eqref{eq:P}.
				\end{enumerate}
			\end{thm}
		\end{mybox}
		\begin{proof}
			\begin{proofitemize}
			\item{thm:KKT:bounded}
				It follows from \cref{thm:decrease} that \(q(x^k)\leq q_{\mu_0}(x^1)<\infty\) holds for every \(k\geq1\).
				Since \(c(x^k)<0\) (because \(x^k\in\dom q_{\mu_{k-1}}\)), one has that \(q(x^k)=q_0(x^k)\), hence that for every \(k\geq1\) \(x^k\) belongs to the sublevel set \(\lev_{\leq q_{\mu_0}(x^1)}q_0\), which is bounded by assumption.

			\item{thm:KKT:feas}
				That \(c(x^\star)\leq0\) follows from \cref{thm:xkfeas} in light of continuity of \(c\).
				Similarly, since \(\seq{q(x^k)}[k\in\N]\) is upper bounded as shown in \cref{thm:decrease}, the inclusion \(x^\star\in\dom q\) owes to lsc of \(q\).

			\item{thm:KKT:ep}
				Among the two possibilities, the algorithm terminates in finite time only if \(\epsilon_{\rm p}=0\) and the returned pair \((x^\star,y^\star)\) satisfies \(\min\set*{-c(x^\star), y^\star}=0\).
				Excluding this ideal situation, we may assume that it runs indefinitely and that consequently \(\mu_k\to0\).
				By \cref{thm:xkfeas,thm:ykfeas}, it is $c(x^k) < 0$ and $y^k \geq 0$ for all $k\in\N$.
				If for some \(\delta>0\) and \(i\in\set{1,\dots,m}\) a subsequence \(\seq{x^k}[k\in K']\) satisfies \(-c_i(x^k)\geq\delta\) for all \(k\in K'\), then \(\seq{b'(c_i(x^k))}[k\in K']\) is bounded and therefore \(y_i^k=\mu_{k-1}b'(c_i(x^k))\to0\) as \(K'\ni k\to\infty\).
				The claim then follows from the arbitrarity of the subsequence.

			\item{thm:KKT:AKKT}
				As shown in assertion \ref{thm:KKT:feas}, \(x^\star\) is feasible.
				Also, \cref{thm:ykfeas,thm:xkopt} together with the fact that $\varepsilon_k \to 0$ ensure that the sequence \(\seq{x^k,y^k}[k\in K]\) satisfies condition \eqref{eq:AKKT:x}.
				Condition \eqref{eq:AKKT:y} follows from assertion \ref{thm:KKT:ep} together with \cref{thm:yk}.

			\item{thm:KKT:KKT}
				Follows from the previous assertion together with \cref{thm:ykbounded}.
			\qedhere
			\end{proofitemize}
		\end{proof}

		\begin{mybox}
			\begin{cor}[Finite termination of \cref{alg:IP}]\label{thm:IP:return}%
				For any strictly feasible starting point \(x^0\) and primal-dual tolerance parameters \(\epsilon_{\rm p},\epsilon_{\rm d}>0\), in finitely many steps \cref{alg:IP} returns an \((\epsilon_{\rm p},\epsilon_{\rm d})\)-KKT optimal point \(x^\star\) for \eqref{eq:P} satisfying \(q(x^\star)\leq q(x^0)\).
			\end{cor}
		\end{mybox}

		Notice that the coercivity assumption of \(q_0\) in \cref{thm:KKT:bounded} needed to ensure boundedness of the sequence generated by \cref{alg:IP} also guarantees that the cost \(q_{\mu_k}\) in each subproblem is level bounded, which is a trivial consequence of the fact that \(q_0\leq q_\mu\) for any \(\mu>0\).
		This in particular guarantees that each subproblem \eqref{eq:Pmu}, for any \(\mu>0\), admits global minimizers.
		Nevertheless, the successful termination of each call to \ipfb{} at \cref{state:IP:x} is independent of whether or not this assumption is met, as demonstrated in \cref{thm:FB:return}, nor is the termination of \cref{alg:IP} affected (as long as strictly positive tolerances \(\epsilon_{\rm p},\epsilon_{\rm d}\) are chosen), as commented in the previous corollary.

	\section{Numerical examples}
		In this section we present some experimental results on an ill-conditioned toy problem to illustrate the numerical behavior of \cref{alg:IP,alg:FB}.
		Then, considering a data analysis task, we investigate the influence of hyperparameters and discuss the performance on larger scale problems.

		To graphically summarize our numerical results and compare different solvers, we display \emph{epi-profiles}, \emph{data profiles}, and \emph{(extended) performance profiles}.
		For $\mathcal{P}$ the set of problems and $\mathcal{S}$ the set of solvers, let $t_{s,p}$ denote the user-defined metric for the computational effort required by solver $s \in \mathcal{S}$ to solve instance $p \in \mathcal{P}$ (lower is better).
		We will monitor the (total) number of gradient evaluations,
		so that the computational overhead triggered by backtracking is fairly accounted for.
		\begin{itemize}
		\item
			\emph{Epi-profiles} display the evaluation metric for individual problems in the problem set $\mathcal{P}$, ordered in such a way that for a user-specified base solver $s\in\mathcal{S}$ the evaluation metric monotonically increases with the problem number.
			The lowest point in each column corresponds to the best solver on the respective instance.

		\item
			\emph{Data profiles} display the cumulative distribution function $f_s \colon [0,\infty) \mapsto [0,1]$ of the evaluation metric, namely
			\[
				f_s(t)
			\coloneqq
				\frac{|\set{p\in\mathcal{P}}[t_{s,p} \leq t]|}{|\mathcal{P}|}
			.
			\]
			Each data profile reports the fraction of problems $f_s(t)$ solved by solver $s$ with a budget $t$ of evaluation metric \cite{more2009benchmarking}, and therefore it is independent of the other solvers.

		\item
			\emph{Extended performance profiles} address the relative performance of solvers \cite[\S4.1]{mahajan2012solving}.
			Let $\tau_{s,p}$ denote the (extended) \emph{performance ratio} of solver $s \in \mathcal{S}$ on a certain instance $p \in \mathcal{P}$ in comparison to the best solver, other than $s$ itself, on that same instance.
			Then, an extended performance profile $\rho_s \colon [0,\infty) \mapsto [0,1]$ is the cumulative distribution function of the performance ratio of solver $s$, namely
			\[
				\rho_s(\tau)
			\coloneqq
				\frac{|\set{p\in\mathcal{P}}[\tau_{s,p} \leq \tau]|}{|\mathcal{P}|}
			\qquad\text{where}\qquad
				\tau_{s,p}
			\coloneqq
				\frac{t_{s,p}}{\min\set{ t_{i,p} }[ i\in\mathcal{S}, i\ne s]}
			.
			\]
			Thus, an extended performance profile indicates the probability (or fraction of problems) $\rho_s(\tau)$ that a given solver $s\in\mathcal{S}$ is faster or slower than any other solver by a given factor $\tau$.
		\end{itemize}%

		\paragraph*{Implementation details}
			We describe here details pertinent to the implementation of \Cref{alg:IP,alg:FB}, defining particular choices left equivocal there,
			such as the initialization and update of algorithmic parameters.
			These numerical features tend to improve the practical performances, without compromising the convergence guarantees established in previous sections.
			\begin{itemize}
			\item
				The initial tolerance \(\varepsilon_0\) for \cref{alg:IP} is chosen adaptively,
				based on the starting point \(x^0\) and barrier parameter \(\mu_0\):
				we set $\varepsilon_0 = \max\set{\epsilon_{\rm d},\kappa_\varepsilon \eta_0}$, where $\kappa_\varepsilon\in(0,1)$ is a user-specified parameter and $\eta_0$ is the norm evaluated for $j=0$ at \cref{state:FB:exit} of \cref{alg:FB} invoked at $(x^0,\mu_0)$.

			\item
				We relax the barrier parameter update rule at \cref{state:IP:epsmu}:
				we set $\mu_{k+1} \gets \mu_k$
				if \((x^{k+1}, y^{k+1})\) satisfies approximate complementarity,
				namely \(\bigl\|\min\set*{ -c(x^{k+1}), y^{k+1} }\bigr\|_\infty \leq \epsilon_{\rm p}\),
				otherwise we reduce the barrier parameter as indicated.

			\item
				The initial stepsize $\gamma_0\in(0,\gamma_g)$ in \cref{alg:FB} is selected adaptively, based on an estimate $L_{\x}$ of $\lip \nabla f_\mu(\x)$.
				We set $\gamma_0 = \alpha / L_{\x}$, where
				$L_{\x} \coloneqq \frac{\|\nabla f_\mu(\x^+) - \nabla f_\mu(\x)\|}{\|\x^+ - \x\|}$
				is a lower bound on the smoothness constant around $\x$.
				The point $\x^+ \coloneqq \x + h$ is obtained by backtracking,
				starting from $h = 1$ and reducing $h$ by a factor $\beta$ until
				$\x^+\in\dom f_\mu$.
				This procedure is well defined since $\x\in\dom f_\mu$ and $c$ is continuous.\footnote{%
					In case $\gamma_g$ is finite the value should be then projected onto $[\delta,\gamma_g-\delta]$ for some $\delta>0$.
					If $L_{\x} = 0$, the choice of \(\gamma_0\) can be arbitrary.
					These minor technicalities are not part of the implementation.
				}

			\item%
				The algorithmic parameters have been set with the following (default) values:
				\(\kappa_\varepsilon = 10^{-2}\),
				\(\mu_0 = 1\),
				\(\theta_\varepsilon = \theta_\mu = \nicefrac{1}{4}\)
				in \cref{alg:IP},
				\(\alpha = 0.9\), \(\beta=\nicefrac{1}{2}\), \(r=1.1\) in \cref{alg:FB}.

			\item%
				At \cref{state:IP:epsmu} of \cref{alg:IP} we always select the respective upper bounds, namely we set $\varepsilon_{k+1} \gets \max\set{\epsilon_{\rm d}, \theta_\varepsilon \varepsilon_k}$ and $\mu_{k+1} \gets \theta_\mu \mu_k$ (or \(\mu_{k+1}=\mu_k\) as described above).

			\item%
				Finally, for constructing the subproblems \eqref{eq:Pmu}, we consider the barrier function $b$ defined by \(b(t)=-\nicefrac{1}{t}\) for \(t<0\), and \(\infty\) otherwise.
				This choice complies with our requirements for a barrier function, having \(b'(t)=\nicefrac{1}{t^2}>0\) for \(t<0\) and \(b\geq b(-\infty)=0\).
			\end{itemize}

			To ensure the reproducibility of the numerical results presented in this paper, our implementation adheres to the steps detailed in \cref{alg:IP,alg:FB},
			incorporating the practical mechanisms just delineated, but 		without introducing any safeguards such as tolerances to mitigate the effects of machine precision.
			Furthermore, the source code of our implementation has been made available on Zenodo at \href{https://doi.org/\TheCodeZenodoDOI}{\textsc{doi}: \TheCodeZenodoDOI}.

		\subsection{Nonsmooth Rosenbrock with inequalities}
			As an illustrative toy example, we consider a two-dimensional optimization problem involving a nonsmooth Rosenbrock-like objective function and inequality constraints.
			Considering the $\ell_p$-quasinorm $\|\cdot\|_p$ with $p \coloneqq \nicefrac{1}{2}$ and a circle with radius $r_C \coloneqq \nicefrac{1}{2}$ centered at $x_C \coloneqq (-\nicefrac{1}{4}, \nicefrac{1}{4})$, it reads
			\begin{equation}\label{eq:rosenbrock}
				\minimize_{x\in\R^2}~
				100 \bigl(x_2 + 1 - (x_1 + 1)^2\bigr)^2 + \|x\|_p^p
				\qquad
				\stt~
				\|x - x_C\|^2 \geq r_C^2 .
			\end{equation}

			\begin{figure}[t]
				\begin{minipage}{0.5\linewidth}
					\vspace*{0pt}%
					\includetikz[width=\linewidth]{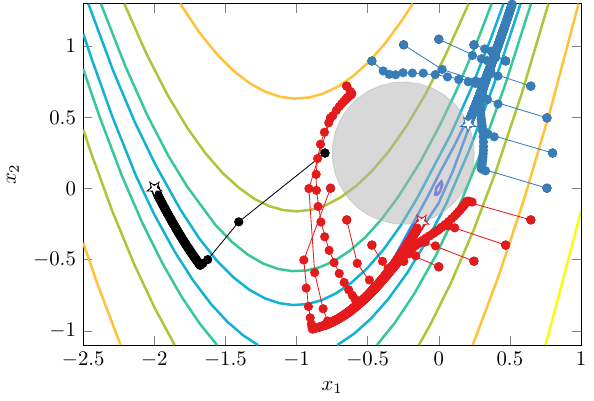}%
					\vspace*{0pt}%
				\end{minipage}
				\hfill
				\begin{minipage}{0.48\linewidth}
					\vspace*{0pt}%
					\caption{%
						Rosenbrock problem \eqref{eq:rosenbrock}: contour lines of the objective function, circular infeasible set (gray), trajectories of inner and outer iterations for different starting points, and limit points thereof (stars).
						Trajectories are colored based on the limit point: $x^{[1]}$ (red), $x^{[2]}$ (blue) or $x^{[3]}$ (black).%
					}%
					\label{fig:rosenbrock_iterates}%
					\vspace*{0pt}%
				\end{minipage}
			\end{figure}

			The proximal mapping of $\|\cdot\|_p^p \colon x \mapsto \sum_{i=1}^2 |x_i|^p$ can be evaluated elementwise based on explicit formulas given in \cite{xu2012regularization,chen2016computing}, namely
			\[
				\Bigl[
					\prox_{\gamma\|{}\cdot{}\|_{\nicefrac12}^{\nicefrac12}}(x)
				\Bigr]_i
			\ni
				\begin{ifcases}
					\frac23\Bigl(
						1+\cos\Bigl(
							\frac23\arccos\bigl(
								-\frac\gamma4
								\bigl(\frac{3}{|x_i|}\bigr)^{\nicefrac32}
							\bigr)
						\Bigr)
					\Bigr)
				&
					|x_i|>\frac32\gamma^{\nicefrac23}
				\\
					0\otherwise.
				\end{ifcases}
			\]
			Furthermore, casting \eqref{eq:rosenbrock} into the form of \eqref{eq:P}, the problem data functions satisfy the conditions in \cref{ass:basic}.
			In particular, $f$ has locally (and not globally) Lipschitz continuous gradient and $g$ is continuous relative to its domain $\dom g = \R^2$.

			We invoked the proposed algorithm on the same problem instance,
			with $\epsilon_{\rm p} = \epsilon_{\rm d} = 10^{-5}$, 			starting from $20$ different (strictly feasible) points $x^0 \in \R^2$.
			These have been generated as $x^0 = (0,\nicefrac{1}{4}) + \nicefrac{4}{5} (\cos\vartheta,\sin\vartheta)$ where $\vartheta \in \R$ is sampled from a uniform grid over $[0,2 \pi]$.

			\begin{figure}[b]
				\centering
				\includetikz[width=.9\linewidth]{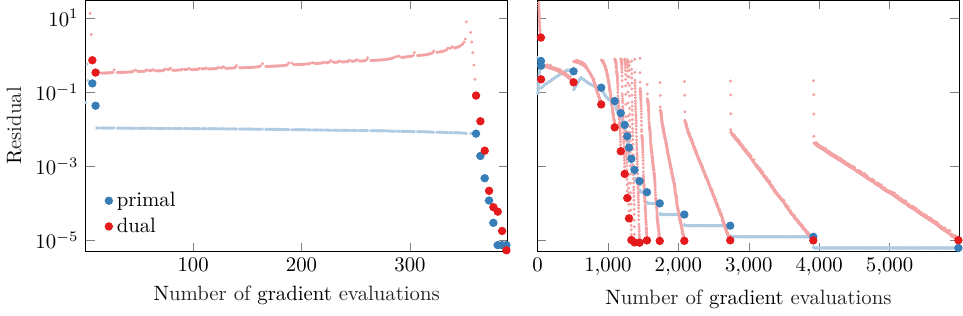}%
				\caption{%
					Rosenbrock problem \eqref{eq:rosenbrock}: comparison of primal and dual residuals against the number of gradient evaluations, for the starting points $x^0 = (-0.8,0.25)$ (left) and $x^0 = (0.8,0.25)$ (right) whose associated limit points are $x^{[3]}$ and $x^{[2]}$, respectively.
					Larger dots correspond to the outer iterations.%
				}%
				\label{fig:rosenbrock_residuals}%
			\end{figure}

			\cref{fig:rosenbrock_iterates,fig:rosenbrock_residuals} summarize the outcomes of these simulations.
			Superimposed to the objective contour lines and the (in)feasible set, the numerical trajectories are depicted in \cref{fig:rosenbrock_iterates}, concatenating over $k=1,2,\ldots$ the iterates $\seq{x^{k,j}}[j\in\N]$ generated by \ipfb.
			Depending on the starting point, \cref{alg:IP} returns one of three stationary points for \eqref{eq:P}, which are indeed the global minimizer $x^{[1]} \approx (-0.12, -0.23)$ or two local minimizers $x^{[2]} \approx (0.21, 0.45)$ and $x^{[3]} \approx (-2.00,0)$, see \cref{fig:rosenbrock_iterates}.
			Notice that the feasible set is not simply connected (hence is nonconvex) and that the constraint is active for two minimizers.
			We stress that the iterates remain strictly feasible while reducing the objective value.

			The algorithm performance in terms of optimality and complementarity measures is illustrated in \cref{fig:rosenbrock_residuals} for two different starting points.
			We monitored the outer dual residual (associated to the inner residual of \cref{state:FB:exit}) and the outer primal residual of \cref{state:IP:exit} at all iterations.
			In accordance with \cref{thm:xkopt}, the dual residual decreases as dictated by the sequence of inner tolerances $\seq{\varepsilon_k}$.
			It is interesting to notice that, even though \cref{thm:KKT:ep} only implies the vanishing of the primal residual, in our simulations it is also monotonically decreasing along outer iterations.

		\subsection{Nonnegative PCA}\label{sec:PCA}%
			Principal component analysis (PCA) aims at estimating the direction of maximal variability of a high-dimensional dataset.
			Arguably the most successful of dimensionality reduction techniques \cite{montanari2016nonnegative},
			classical PCA aims to recover a signal $z$ from finding the eigenvector that corresponds to the largest eigenvalue of a given matrix $Z$ \cite{liu2020simple}.
			A recurring idea is to use additional structural information about the principal eigenvector, such as its signature or sparsity \cite{montanari2016nonnegative}.
			Here we impose nonnegativity of entries as prior knowledge, and solve PCA restricted to the positive orthant:
			\begin{equation}
				\label{eq:nonneg_pca}
				\maximize_{x\in\R^n}~ x^\top Z x \qquad
				\stt~ \| x \| = 1 ,\ x \geq 0 .
			\end{equation}
			This task falls within the scope of \eqref{eq:P}, with $f(x) \coloneqq - x^\top Z x$, $g(x) \coloneqq \indicator_{\|\cdot\| = 1}(x)$, and $c(x) = - x$.
			Nonnegative PCA is an NP-hard nonconvex problem \cite{montanari2016nonnegative} that cannot be addressed by standard SVD.

			\paragraph*{Setup}
				We synthetically generate problem data following \cite{liu2020simple}.
				For a problem size $n\in\N$,
				let $Z = \sqrt{\sigma} z z^\top + N\in\R^{n\times n}$, where $N\in\R^{n\times n}$ is a random symmetric noise matrix and $\sigma>0$ is the signal-to-noise ratio (SNR).
				The off-diagonal entries of $N$ follow a Gaussian distribution $\mathcal{N}(0,\nicefrac{1}{n})$ and its diagonal entries follow a Gaussian distribution $\mathcal{N}(0,\nicefrac{2}{n})$.
				Furthermore, we let the support $S \subseteq \set{1,\ldots,n}$ of the true principal direction $z$ be uniformly random, with cardinality $|S| = \lfloor s n \rfloor$, and
				set $z_i = \nicefrac{1}{\sqrt{|S|}}$ if $i\in S$, $z_i=0$ otherwise.
				We consider some dimensions $n$ and,
				for each dimension,
				the set of problems parametrized by $\sigma \in \set{0.05, 0.1, 0.25, 0.5, 1.0}$ and $s \in \set{0.1, 0.3, 0.7, 0.9}$, which control the noise and sparsity level, respectively.
				A strictly feasible starting point $x_0$ is generated by sampling a uniform distribution over $[0,3]^n$ and projecting onto $\dom g = \set{x\in\R^n}[\|x\| = 1]$.
				There are 5 choices for $\sigma$, 4 for $s$, and, for each set of parameters, 5 instances are generated with different problem data $Z$ and starting point $x_0$.
				Overall, each solver-settings pair is invoked on 100 different instances for each dimension $n$.

			\begin{figure}[t]
				\centering
				\includetikz[width=.9\linewidth]{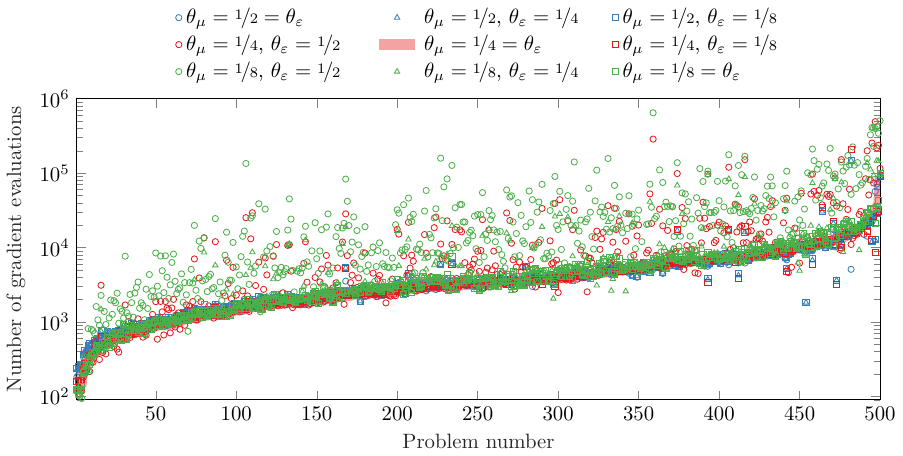}%
				\caption{%
					Nonnegative PCA problem \eqref{eq:nonneg_pca}:
					comparison for different barrier parameter and inner tolerance reduction factors $\theta_\mu,\theta_\varepsilon\in(0,1)$.
					Epi-profiles ordered in relation to the default values $\theta_\mu = \nicefrac14 = \theta_\varepsilon$ (red thick line).%
				}%
				\label{fig:nonneg_pca_thetas_epiprof}%
			\end{figure}

			\paragraph*{Hyperparameters tuning}
				\Cref{alg:IP,alg:FB} are controlled by several hyperparameters, such as the initial barrier parameter $\mu_0$, reduction factors $\theta_\mu, \theta_\varepsilon$, and the regret factor $r$.
				Investigating the influence of hyperparameters is not only interesting to effectively tune the solvers, but also to appreciate how sensitive (or robust) the performance is with respect to their values.

			\begin{figure}[t]
				\centering
				\includetikz[width=.9\linewidth]{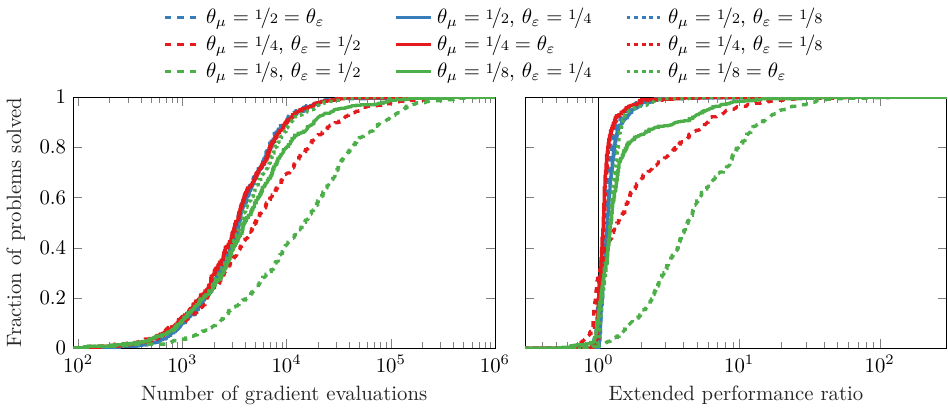}%
				\caption{%
					Nonnegative PCA problem \eqref{eq:nonneg_pca}:
					comparison for different barrier parameter and inner tolerance reduction factors $\theta_\mu,\theta_\varepsilon\in(0,1)$.
					Data profiles (top) and extended performance profiles (bottom) relative to the number of gradient evaluations.%
				}%
				\label{fig:nonneg_pca_thetas_timeprof_perfprof}%
			\end{figure}

			We now focus on the effect of $\theta_\mu,\theta_\varepsilon\in(0,1)$,
			considering problem dimensions $n \in \set{10, 15, 20, 25, 30}$ and all combinations of $\theta_\mu,\theta_\varepsilon\in\set{\nicefrac12, \nicefrac14, \nicefrac18}$, for a total of 4500 calls to \cref{alg:IP},
			with tolerances $\epsilon_{\rm p} = \epsilon_{\rm d} = 10^{-3}$.
			Lower values of $\theta_\mu$ ($\theta_\varepsilon$) yield a faster decrease of the barrier parameters $\mu_k$ (inner tolerances $\varepsilon_k$) toward zero.

			All instances are solved up to the desired primal-dual tolerances.
			The results are graphically summarized in \cref{fig:nonneg_pca_thetas_epiprof,fig:nonneg_pca_thetas_timeprof_perfprof},
			showing that the majority of selected tunings yield comparable results.
			The settings $(\theta_\mu,\theta_\varepsilon) = (\nicefrac18,\nicefrac14)$, $(\nicefrac14,\nicefrac12)$, and $(\nicefrac18,\nicefrac12)$ are increasingly worse,
			whereas $(\theta_\mu,\theta_\varepsilon) = (\nicefrac14,\nicefrac14)$ seems to dominate.
			This value agrees with the default settings chosen for the solver, as mentioned in the beginning of this section.

			Let us now examine the influence of the regret factor $r \geq 1$ in \cref{alg:FB},
			considering the same problem setup and the values $r\in\set{1, 1.1, 1.25, 1.5}$,
			for a total of 2000 calls to \cref{alg:IP},
			including the adaptive variant described below.
			As commented in the beginning of \cref{sec:FB}, higher values of $r$ allow the stepsize to recover faster from low values that compromise convergence speed, when the local geometry of \(f\) allows.
			On the other hand, lower values of \(r\) reduce the number of backtrackings at every step, and thus the number of gradient evaluations per iteration.
			However, other than keeping $r$ constant, it is possible to consider any sequence \(\seq{r_j}[j\in\N]\subset[1,\infty)\),
			as mentioned in the discussion after \cref{thm:PG}.
			This motivates testing also \cref{alg:FB} with an adaptive regret:
			on the line of \cite{malitsky2020adaptive},
			we consider the sequence generated by $r_j=\sqrt{1+\nicefrac{\gamma_{j-1}}{\gamma_{j-2}}}$ for all $j\geq 2$,
			with the initialization $\gamma_j=r_j\gamma_{j-1}$ at \cref{state:FB:init}.

\begin{figure}[t]
				\centering
				\includetikz[width=.9\linewidth]{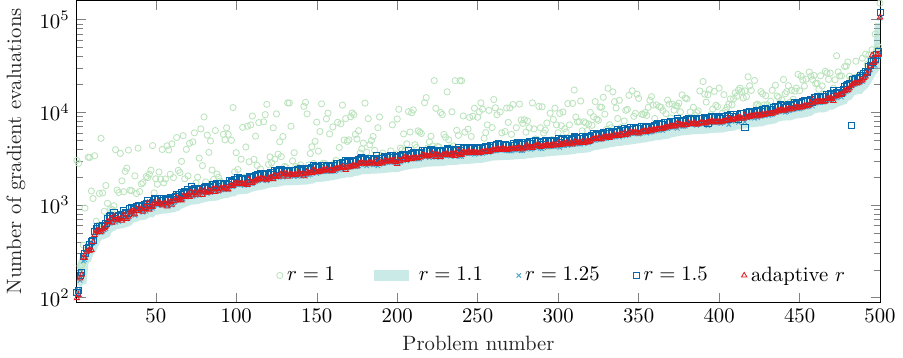}%
				\caption{%
					Nonnegative PCA problem \eqref{eq:nonneg_pca}:
					comparison for different regret factors $r\geq 1$.
					Epi-profiles ordered in relation to the default value $r = 1.1$ (thick line).%
				}%
				\label{fig:nonneg_pca_regret_epiprof}%
			\end{figure}

\begin{figure}[bh]
				\centering
				\includetikz[width=.9\linewidth]{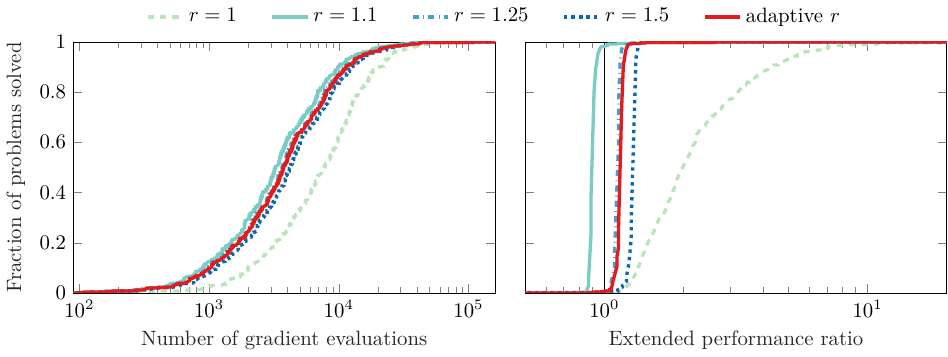}%
				\caption{%
						Nonnegative PCA problem \eqref{eq:nonneg_pca}:
						comparison for different regret factors $r\geq 1$.
						Data profiles (top) and extended performance profiles (bottom) relative to the number of gradient evaluations.%
				}%
				\label{fig:nonneg_pca_regret_timeprof_perfprof}%
			\end{figure}

			All instances are solved up to the desired primal-dual tolerances
			and computational results are graphically summarized in \cref{fig:nonneg_pca_regret_epiprof,fig:nonneg_pca_regret_timeprof_perfprof}.
			According to these profiles, a suitable tuning for the regret factor in \cref{alg:FB} appears to be around the value $r = 1.1$,
			in agreement with the default settings chosen for the solver.
			These results illustrate the significant potential benefits of a regret factor \(r>1\),
			as revealed by the considerable gap with the monotone stepsize initialization ($r=1$).
			Furthermore, all the tested values of \(r>1\) yield consistent improvements over the choice \(r=1\), indicating that the good performance of \cref{alg:FB} may be robust with respect to the regret factor for values of $r$ strictly larger than (but close to) 1.
			This is also true for the adaptive choice, suggesting that it could also constitute a conveniently parameter-free strategy.

			\begin{figure}[t]
				\centering
				\includetikz[width=.9\linewidth]{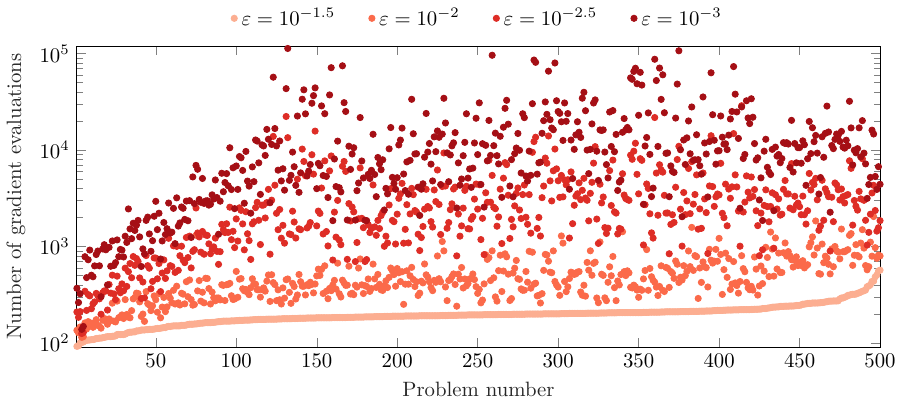}%
				\caption{%
					Nonnegative PCA problem \eqref{eq:nonneg_pca}:
					comparison for increasing accuracy requirements (decreasing tolerances $\epsilon_{\rm p} = \epsilon_{\rm d} = \varepsilon$).
					Epi-profiles ordered in relation to $\varepsilon=10^{-1.5}$.%
				}%
				\label{fig:nonneg_pca_epiprof}%
			\end{figure}
			\begin{figure}[t]
				\centering
				\includetikz[width=.9\linewidth]{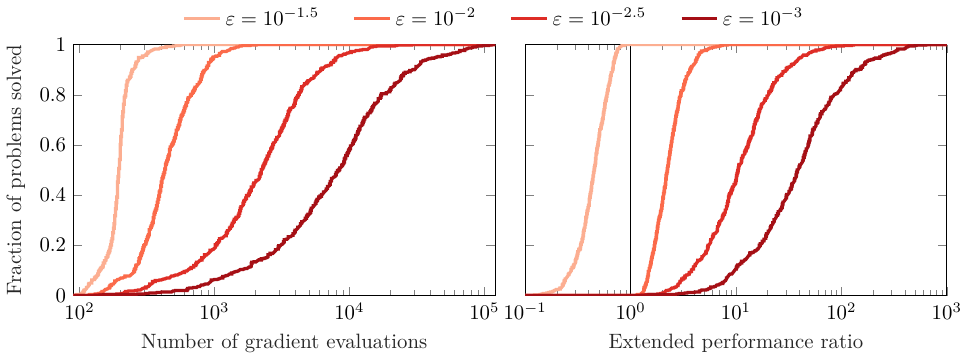}%
				\caption{%
					Nonnegative PCA problem \eqref{eq:nonneg_pca}:
					comparison for increasing accuracy requirements (decreasing tolerances $\epsilon_{\rm p} = \epsilon_{\rm d} = \varepsilon$).
					Data profiles (left) and extended performance profiles (right) relative to the number of gradient evaluations.%
				}%
				\label{fig:nonneg_pca_timeprof_perfprof}%
			\end{figure}

			\paragraph*{Problem size and tolerance}
			To investigate scalability and influence of accuracy requirements,
			we consider instances of \eqref{eq:nonneg_pca} with dimensions $n \in \set{10, \lceil 10^{1.5} \rceil, 10^2, \lceil 10^{2.5} \rceil, 10^3}$ and tolerances $\epsilon_{\rm p} = \epsilon_{\rm d} = \varepsilon \in \set{10^{-1.5}, 10^{-2}, 10^{-2.5}, 10^{-3}}$.
			Each of these tolerance parameters is tested on 500 problem instances, for a total of 2000 calls to \cref{alg:IP}.

\begin{figure}[t]
				\centering
				\includetikz[width=.9\linewidth]{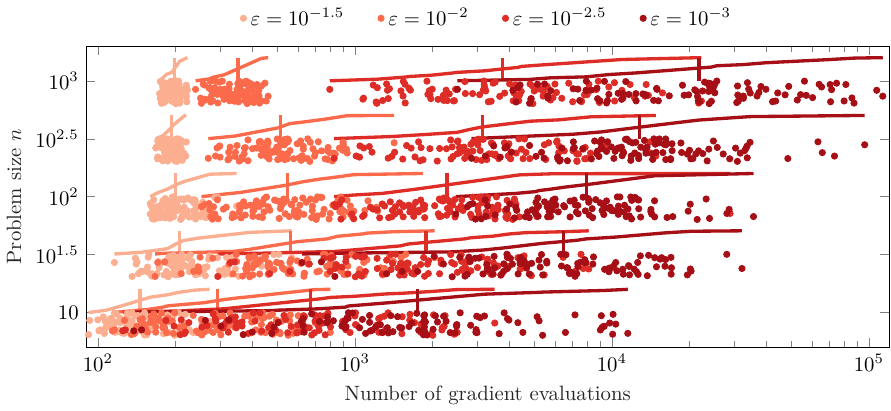}%
				\caption{%
					Nonnegative PCA problem \eqref{eq:nonneg_pca}:
					comparison for increasing accuracy requirements (decreasing tolerances $\epsilon_{\rm p} = \epsilon_{\rm d} = \varepsilon$) and problem sizes $n$.
					Combination of jitter plot (dots) and cumulative distribution function estimate (solid line) with median value (vertical line).%
				}%
				\label{fig:nonneg_pca_jittercdf}%
			\end{figure}

			All instances are solved up to the desired primal-dual tolerances.
			The results are graphically summarized in \cref{fig:nonneg_pca_epiprof,fig:nonneg_pca_timeprof_perfprof},
			where it is clear that stricter tolerances demand more effort, as expected.
			However, it is interesting to look at how the computational cost significantly increases with the accuracy requirement, because of the slow tail convergence typical of first-order methods such as \ipfb.
			The influence of tolerance and problem size is depicted in \cref{fig:nonneg_pca_jittercdf},
			which displays for each pair $(n,\varepsilon)$ the number of gradient evaluations with a jitter plot
			and
			reports an estimate of the cumulative distribution function with the associated median value.\footnote{%
				Jitter plots offer a simple way of visualizing the distribution of numerical values over categories.
				Sample values are plotted as dots along one axis, shifted randomly along the other axis;
				the jittering has no meaning in itself data-wise, but allows a better view of overlapping data points.
				Jitter plots are complemented with the cumulative distribution function, as opposed to the probability density function, since a robust estimate of the former does not require additional assumptions.
				The combined plot thus conveys information on the number of data points and their density distribution in an honest and comprehensible format.%
			}
			This chart visualizes how problem size and accuracy requirement affect the solution process,
			and reveals the stark effect of both $n$ and $\varepsilon$.

	\section{Conclusions}
		We proposed an interior point (IP) method for nonsmooth minimization subject to smooth inequality constraints, where the inner barrier subproblems are addressed by means of proximal gradient iterations.
		The methodology is an extension to a fully nonconvex setting of the PIPA algorithm proposed in \cite{chouzenoux2020proximal}, and aims at bridging the gap between IP and proximal algorithms, the former being the methods of choice for coping with complex constraints and the latter being well suited for large-scale nonsmooth problems.
		The result is a warm-startable iterative scheme whose output are approximate KKT-optimal pairs for the problem.
		Our analysis of proximal gradient iterations is novel, offering weaker conditions to ensure convergence results in the fully nonconvex setting.

		Despite the benefits of adopting nonmomontone stepsize sequences demonstrated by our numerical simulations , the method suffers from the slow tail convergence that is typical of first-order methods.
		These observations motivate future research directions toward integrating the methodology with more adaptive and higher-order schemes.
		While the direct adoption of accelerated solvers along the lines of \cite{themelis2018forward,demarchi2022proximal} seems far from trivial, variable-metric or proximal-Newton approaches could be viable options for coping with the ill-conditioning inherent to the barrier subproblems, as observed in \cite{chouzenoux2020proximal}.
		Other interesting developments include gaining a deeper understanding on the choice of barrier parameters and inner tolerances to improve convergence and output quality.
		Finally, a non-asymptotic analysis of \cref{alg:IP,alg:FB} is left for future work, to shed light on whether there is a uniform upper bound on the number of steps, or under which conditions.
		In particular, as \cref{cond:FB:domfmu} affects the linesearch procedure, maintaining strict feasibility seems to hinder complexity estimates in the nonconvex setting of \cref{ass:basic}, suggesting that additional assumptions may be required for the purpose.

	\phantomsection
	\addcontentsline{toc}{section}{References}%
	\small
	\bibliographystyle{plain}
	\bibliography{ipprox.bib}

\end{document}